\documentclass{article}
\usepackage[T1]{fontenc}
\usepackage{lmodern}
\usepackage{graphicx, amsmath, amsfonts, amsthm, amssymb, tikz}
\usepackage[nocompress]{cite}
\usepackage[hidelinks]{hyperref}
\usepackage{parskip}

\title{Combinatorial Zeta Functions of Type $\widetilde A_2$}
\author{Ming-Hsuan Kang and  Chian-Jen Wang}
\date{August 2026}
\usetikzlibrary{calc,arrows.meta}
\newtheoremstyle{plainsp}{0.7\baselineskip plus 0.2\baselineskip minus 0.1\baselineskip}%
  {0.7\baselineskip plus 0.2\baselineskip minus 0.1\baselineskip}%
  {\itshape}{}{\bfseries}{.}{.5em}{}
\newtheoremstyle{definitionsp}{0.7\baselineskip plus 0.2\baselineskip minus 0.1\baselineskip}%
  {0.7\baselineskip plus 0.2\baselineskip minus 0.1\baselineskip}%
  {\normalfont}{}{\bfseries}{.}{.5em}{}
\newtheoremstyle{remarksp}{0.7\baselineskip plus 0.2\baselineskip minus 0.1\baselineskip}%
  {0.7\baselineskip plus 0.2\baselineskip minus 0.1\baselineskip}%
  {\normalfont}{}{\itshape}{.}{.5em}{}

\theoremstyle{plainsp}
\newtheorem{theorem}{Theorem}[section]

\newtheorem{proposition}[theorem]{Proposition}
\newtheorem{lemma}[theorem]{Lemma}

\theoremstyle{definitionsp}
\newtheorem{definition}[theorem]{Definition}
\newtheorem{example}[theorem]{Example}

\theoremstyle{remarksp}
\newtheorem{remark}[theorem]{Remark}

\newcommand{\Id}{\mathrm{Id}}

\newcommand{\DrawThetaLink}{%
  \coordinate (L)  at (0,0);
  \coordinate (R)  at (4,0);

  \coordinate (T1) at (1.3,0.85);
  \coordinate (T2) at (2.7,0.85);

  \coordinate (M1) at (1.3,0);
  \coordinate (M2) at (2.7,0);

  \coordinate (B1) at (1.3,-0.85);
  \coordinate (B2) at (2.7,-0.85);

  \draw[thick] (L) -- (T1) -- (T2) -- (R);
  \draw[thick] (L) -- (M1) -- (M2) -- (R);
  \draw[thick] (L) -- (B1) -- (B2) -- (R);

  \foreach \v in {L,T2,M2,B2} {
    \filldraw[fill=white,draw=black] (\v) circle (2.2pt);
  }

  \foreach \v in {R,T1,M1,B1} {
    \filldraw[fill=black,draw=black] (\v) circle (2.2pt);
  }
}

\newcommand{\DrawHexLink}{%
  \foreach \k in {0,...,5} {
    \coordinate (H\k) at ({cos(60*\k)},{sin(60*\k)});
  }

  \foreach \a/\b in {0/1,1/2,2/3,3/4,4/5,5/0} {
    \draw[thick] (H\a) -- (H\b);
  }

  \foreach \k in {0,2,4} {
    \filldraw[fill=white,draw=black] (H\k) circle (2.2pt);
  }

  \foreach \k in {1,3,5} {
    \filldraw[fill=black,draw=black] (H\k) circle (2.2pt);
  }
}

\newcommand{\DrawCorrectTorusCopy}[4]{%
\begin{scope}[shift={#1}, x={(0.75cm,0cm)}, y={(-0.375cm,0.649519cm)},
              line join=round, line cap=round]

  \foreach \i in {0,...,3} {
    \foreach \j in {0,...,3} {
      \coordinate (#3P\i\j) at (\i,\j);
    }
  }

  \fill[#2!8] (#3P00) -- (#3P30) -- (#3P33) -- (#3P03) -- cycle;
  \fill[#4] (#3P10) -- (#3P20) -- (#3P23) -- (#3P13) -- cycle;

  \foreach \j in {0,...,3} {
    \draw[gray!55] (#3P0\j) -- (#3P3\j);
  }
  \foreach \i in {0,...,3} {
    \draw[gray!55] (#3P\i0) -- (#3P\i3);
  }

  \foreach \i in {0,...,2} {
    \foreach \j in {0,...,2} {
      \pgfmathtruncatemacro{\ip}{\i+1}
      \pgfmathtruncatemacro{\jp}{\j+1}
      \draw[gray!55] (#3P\i\j) -- (#3P\ip\jp);
    }
  }

  \draw[very thick] (#3P00) -- (#3P30) -- (#3P33) -- (#3P03) -- cycle;

  \foreach \i in {0,...,3} {
    \foreach \j in {0,...,3} {
      \pgfmathtruncatemacro{\c}{mod(\i+\j,3)}
      \ifnum\c=0
        \fill[red!75] (#3P\i\j) circle (2.2pt);
      \fi
      \ifnum\c=1
        \fill[blue!75] (#3P\i\j) circle (2.2pt);
      \fi
      \ifnum\c=2
        \fill[green!65!black] (#3P\i\j) circle (2.2pt);
      \fi
    }
  }

  \draw[->, very thick] (0,-0.45) -- (3,-0.45);
  \draw[->, very thick] (0,3.45) -- (3,3.45);
  \draw[->, very thick] (-0.45,0) -- (-0.45,3);
  \draw[->, very thick] (3.45,0) -- (3.45,3);
\end{scope}
}

\newcommand{\DrawKleinBottleCopy}[3]{%
\begin{scope}[shift={#1}, x={(0.58cm,0cm)}, y={(-0.29cm,0.5023cm)},
              line join=round, line cap=round]

\fill[#2!4] (-2.5,0) -- (1.5,2) -- (1.5,11) -- (-2.5,9) -- cycle;

\fill[#3] (-0.5,1) -- (1.5,2) -- (1.5,11) -- (-0.5,10) -- cycle;

\begin{scope}
\clip (-2.5,0) -- (1.5,2) -- (1.5,11) -- (-2.5,9) -- cycle;

\foreach \j in {-1,...,12} {
  \draw[gray!55] (-5,\j) -- (5,\j);
}

\foreach \i in {-5,...,4} {
  \draw[gray!55] (\i,-1) -- (\i,12);
}

\foreach \k in {-12,...,8} {
  \draw[gray!55] (\k,-1) -- ({\k+13},12);
}

\foreach \i in {-3,...,2} {
  \foreach \j in {0,...,11} {
    \pgfmathtruncatemacro{\c}{mod(\i+\j+300,3)}
    \ifnum\c=0
      \fill[red!75] (\i,\j) circle (2.1pt);
    \fi
    \ifnum\c=1
      \fill[blue!75] (\i,\j) circle (2.1pt);
    \fi
    \ifnum\c=2
      \fill[green!65!black] (\i,\j) circle (2.1pt);
    \fi
  }
}
\end{scope}

\draw[very thick]
  (-2.5,0) -- (1.5,2) -- (1.5,11) -- (-2.5,9) -- cycle;

\draw[dashed, very thick]
  (-0.5,1) -- (-0.5,10);

\draw[->, very thick]
  (-2.5,-0.45) -- (1.5,1.55);
\node[font=\small] at (-0.45,-0.55) {$T=(4,2)$};

\draw[->, very thick]
  (-3,0) -- (-3,9);
\node[font=\small, rotate=90] at (-4.25,4.9) {$g^2=(0,9)$};

\draw[->, very thick]
  (-1.65,3.2) to[out=20,in=210] (0.65,8.85);
\node[font=\small] at (-0.20,6.25) {$g$};

\end{scope}
}

\newcommand{\DrawAffineApartmentPatch}{%
\begin{scope}[x={(0.58cm,0cm)}, y={(-0.29cm,0.5023cm)},
              line join=round, line cap=round]
  \fill[blue!3] (-0.5,-0.5) -- (5.5,-0.5) -- (5.5,4.5) -- (-0.5,4.5) -- cycle;
  \begin{scope}
    \clip (-0.5,-0.5) -- (5.5,-0.5) -- (5.5,4.5) -- (-0.5,4.5) -- cycle;
    \foreach \j in {-1,...,5} { \draw[gray!55] (-1,\j) -- (6,\j); }
    \foreach \i in {-1,...,6} { \draw[gray!55] (\i,-1) -- (\i,5); }
    \foreach \k in {-6,...,6} { \draw[gray!55] (\k,-1) -- ({\k+6},5); }
    \foreach \i in {0,...,5} {
      \foreach \j in {0,...,4} {
        \pgfmathtruncatemacro{\c}{mod(\i+\j,3)}
        \ifnum\c=0
          \fill[red!75] (\i,\j) circle (2.1pt);
        \fi
        \ifnum\c=1
          \fill[blue!75] (\i,\j) circle (2.1pt);
        \fi
        \ifnum\c=2
          \fill[green!65!black] (\i,\j) circle (2.1pt);
        \fi
      }
    }
  \end{scope}
  \draw[very thick] (-0.5,-0.5) -- (5.5,-0.5) -- (5.5,4.5) -- (-0.5,4.5) -- cycle;
  \draw[->, thick] (0,0) -- (1,0);
  \node[font=\scriptsize, above] at (0.9,-.1) {$(1,0)$};
  \draw[->, thick] (0,0) -- (0,1);
  \node[font=\scriptsize, above] at (-.2,0.9) {$(0,1)$};
  \draw[->, thick] (0,0) -- (1,1);
  \node[font=\scriptsize, above] at (1.1,0.9) {$(1,1)$};
\end{scope}%
}

\newcommand{\DrawTorusNineStripBlock}[3]{%
\begin{scope}[shift={#1}, x={(0.46cm,0cm)}, y={(-0.23cm,0.3984cm)},
              line join=round, line cap=round]
  \fill[#2!8] (0,0) -- (3,0) -- (3,9) -- (0,9) -- cycle;
  \fill[#3] (1,0) -- (2,0) -- (2,9) -- (1,9) -- cycle;
  \foreach \j in {0,...,9} { \draw[gray!55] (0,\j) -- (3,\j); }
  \foreach \i in {0,...,3} { \draw[gray!55] (\i,0) -- (\i,9); }
  \foreach \i in {0,...,2} {
    \foreach \j in {0,...,8} {
      \pgfmathtruncatemacro{\ip}{\i+1}
      \pgfmathtruncatemacro{\jp}{\j+1}
      \draw[gray!55] (\i,\j) -- (\ip,\jp);
    }
  }
  \draw[very thick] (0,0) -- (3,0) -- (3,9) -- (0,9) -- cycle;
  \foreach \i in {0,...,3} {
    \foreach \j in {0,...,9} {
      \pgfmathtruncatemacro{\c}{mod(\i+\j,3)}
      \ifnum\c=0
        \fill[red!75] (\i,\j) circle (1.9pt);
      \fi
      \ifnum\c=1
        \fill[blue!75] (\i,\j) circle (1.9pt);
      \fi
      \ifnum\c=2
        \fill[green!65!black] (\i,\j) circle (1.9pt);
      \fi
    }
  }
  \draw[->, thick] (0,-0.55) -- (3,-0.55);
  \draw[->, thick] (0,9.55) -- (3,9.55);
  \draw[->, thick] (-0.55,0) -- (-0.55,9);
  \draw[->, thick] (3.55,0) -- (3.55,9);
\end{scope}%
}

\newcommand{\DrawKleinDoubleCoverWithLiftedStrips}[3]{%
\begin{scope}[shift={#1}, x={(0.42cm,0cm)}, y={(-0.21cm,0.3637cm)},
              line join=round, line cap=round]
  \coordinate (A) at (-2.5,0);
  \coordinate (B) at (1.5,2);
  \coordinate (C) at (1.5,11);
  \coordinate (D) at (-2.5,9);
  \fill[#2!8] (A) -- (B) -- (C) -- (D) -- cycle;
  \begin{scope}
    \clip (A) -- (B) -- (C) -- (D) -- cycle;
    \fill[#3] (0,-1) -- (1,-1) -- (1,12) -- (0,12) -- cycle;
    \fill[#3] (-2,-1) -- (-1,-1) -- (-1,12) -- (-2,12) -- cycle;
    \foreach \j in {-1,...,12} { \draw[gray!55] (-5,\j) -- (5,\j); }
    \foreach \i in {-5,...,4} { \draw[gray!55] (\i,-1) -- (\i,12); }
    \foreach \k in {-12,...,8} { \draw[gray!55] (\k,-1) -- ({\k+13},12); }
    \foreach \i in {-3,...,2} {
      \foreach \j in {0,...,11} {
        \pgfmathtruncatemacro{\c}{mod(\i+\j+300,3)}
        \ifnum\c=0
          \fill[red!75] (\i,\j) circle (1.8pt);
        \fi
        \ifnum\c=1
          \fill[blue!75] (\i,\j) circle (1.8pt);
        \fi
        \ifnum\c=2
          \fill[green!65!black] (\i,\j) circle (1.8pt);
        \fi
      }
    }
  \end{scope}
  \draw[very thick] (A) -- (B) -- (C) -- (D) -- cycle;
  \draw[dashed, thick] (-0.5,1) -- (-0.5,10);
  \draw[->, very thick] (1,3) to[out=160,in=20] (-2,6);
  \node[font=\small] at (-0.55,5.0) {$g$};
  \draw[->, thick] (-2.5,-0.55) -- (1.5,1.45);
  \node[font=\scriptsize] at (-0.4,-0.55) {$T$};
  \draw[->, thick] (-3.1,0) -- (-3.1,9);
  \node[font=\scriptsize, rotate=90] at (-3.65,4.5) {$g^2$};
\end{scope}%
}

\begin{document}

\maketitle

\begin{abstract}
We develop a combinatorial zeta theory for finite colored triangular
complexes $X$ whose vertex links satisfy a local hexagon condition.
Their edge and chamber zeta functions $Z_1(X,u)$ and $Z_2(X,u)$ are
Euler products over primitive closed geodesics. The M\"obius functions
of the augmented local incidence posets define a vertex polynomial
$M(X,u)$. We prove the identity
\[
\frac{Z_1(X,u)}{Z_2(X,-u)}
=
(1-u^3)^{\chi(X)}M(X,u)^{-1},
\]
where $\chi(X)$ is the Euler characteristic of $X$.
The proof reconstructs the cohomological operator framework for finite
$\mathrm{PGL}_3$ quotients from the local link data: the degree operators
recover the three determinant factors and together form a cochain map
homotopic to multiplication by $1-u^3$. Sylvester's identity gives the
resulting alternating determinant relation directly. All constructions
use only the finite incidence structure of $X$; no ambient group, field,
or building is required.
\end{abstract}

\section{Introduction}

The classical Ihara zeta function of a finite undirected graph is defined as an Euler product over primitive closed geodesics. A central feature of the theory is the equivalence between two determinant formulas: the Ihara formula, which expresses the zeta function via the vertex adjacency matrix, and the Hashimoto formula, which encodes the dynamics via the edge-transition matrix \cite{Ihara1966,Hashimoto1989,Bass1992}. This fundamental identity bridges the global spectral data of the graph with its local geometric walk dynamics.

A connected $k$--regular graph is a finite quotient of the $k$--regular tree. When $k-1$ is a prime power, from an arithmetic perspective, this tree can be viewed as the Bruhat--Tits tree associated with $\mathrm{PGL}_2$ over a non-Archimedean local field, representing the rank-one case of building theory. Extending this framework to rank two, such as quotients of Bruhat--Tits buildings associated with $\mathrm{PGL}_3$ or other groups, has attracted significant attention \cite{KangLiWang2010,KangLi2014,KangLiArtinPGL3,DeitmarKang2013,FangLiWang2013}; very recently, this program has also been extended to non-uniform, finite-volume $\mathrm{PGL}_3$ quotients, where an Ihara--Bass type determinant formula for a chamber zeta function was obtained via a chamber transfer operator \cite{HongKwon2025}. In this setting, the analogue of the zeta identity contains three interacting factors: an edge zeta function, a chamber zeta function, and a vertex $L$--factor. The relation among these three functions was established by exploiting the ambient algebraic structure, through Hecke
operators, representation theory, or cohomological determinant
arguments. For higher-dimensional complexes, Kang--Yu \cite{KangYu2026} recently established a zeta identity for finite quotients of the Bruhat--Tits building of $\mathrm{PGL}_n$ using a cohomological method.

A related spectral direction is the theory of Ramanujan complexes.
Lubotzky, Samuels, and Vishne constructed such complexes as finite
quotients of affine buildings \cite{LSV2005Ramanujan,LSV2005Explicit}.
Closer to the dynamics considered here, Lubetzky, Lubotzky, and
Parzanchevski studied geodesic-flow operators on Ramanujan complexes
and established a Riemann hypothesis for their associated zeta functions
\cite{LLP2020}. These results concern the spectral properties of building
quotients; our aim is to establish the determinant identity from local
incidence data, without a Ramanujan assumption.

\begin{figure}[ht]
\centering
\begin{tikzpicture}[
    scale=1.15,
    line join=round,
    line cap=round,
    >=Latex
]
  \coordinate (v2) at (0,0);
  \coordinate (v0) at (-2,0);
  \coordinate (v1) at (-1,1.732);
  \coordinate (v3) at (1,1.732);
  \coordinate (v4) at (2,0);

  \fill[blue!12] (v2) -- (v0) -- (v1) -- cycle;
  \fill[blue!22] (v2) -- (v1) -- (v3) -- cycle;
  \fill[blue!12] (v2) -- (v3) -- (v4) -- cycle;

  \draw[thick] (v0) -- (v1) -- (v3) -- (v4);
  \draw[thick] (v2) -- (v0);
  \draw[thick] (v2) -- (v1);
  \draw[thick] (v2) -- (v3);
  \draw[thick] (v2) -- (v4);

  \draw[very thick,red!70!black,->]
    (v0) -- (v2);
  \draw[very thick,red!70!black,->]
    (v2) -- (v4);

  \foreach \p in {v0,v1,v2,v3,v4} {
    \fill (\p) circle (2pt);
  }

  \node[left] at (v0) {$v_0$};
  \node[above left] at (v1) {$v_1$};
  \node[below] at (v2) {$v_2$};
  \node[above right] at (v3) {$v_3$};
  \node[right] at (v4) {$v_4$};

  \node at (-1.05,0.55) {$F_0$};
  \node at (0,0.95) {$F_1$};
  \node at (1.05,0.55) {$F_2$};

  \draw[->,blue!70!black,thick]
    (-0.72,1.20) to[out=20,in=160] (0.72,1.20);
\end{tikzpicture}
\caption{A three-chamber fan exhibiting edge and chamber successor motion.}
\label{fig:intro-three-chamber-fan}
\end{figure}
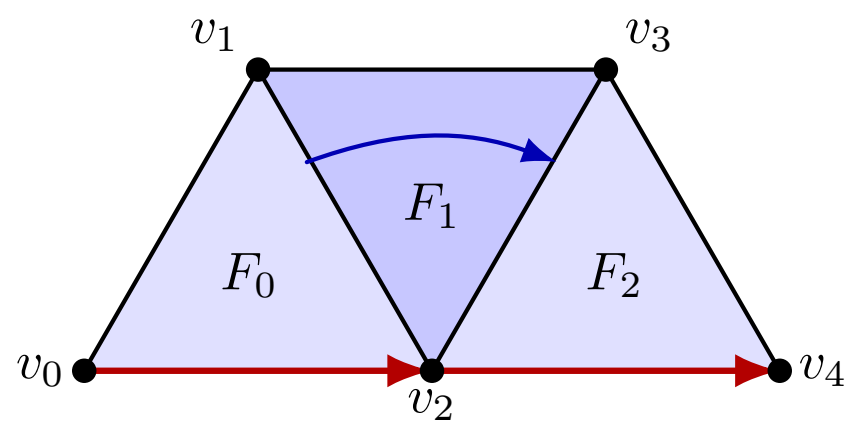

The purpose of the present paper is to isolate the local incidence
mechanism behind the $\mathrm{PGL}_3$ identity, established in \cite{KangLiWang2010, KangLi2014, KangLiArtinPGL3} from different viewpoints, and to formulate it entirely
in finite combinatorial terms, so as to extend the zeta identity to a more general setting explained below.  In a triangular complex, straight motion
is no longer determined by avoiding the reversal of a single edge; the
surrounding chambers record the local direction.  Consider three
consecutive chambers around $v_2$,
\(F_0=\{v_0,v_1,v_2\}\), \(F_1=\{v_1,v_2,v_3\}\), and
\(F_2=\{v_2,v_3,v_4\}\),
as shown in Figure~\ref{fig:intro-three-chamber-fan}.  When these
chambers form one half of a tiled hexagon, the two boundary edges at
$v_2$ correspond to opposite vertices in the \emph{link} $\mathrm{Lk}(v_2)$
--- the graph recording the edges and chambers at $v_2$, defined
formally in Section~\ref{s:local_geometry}.  Consequently,
the directed path continues in the same straight direction through
$v_2$, and we regard \((v_0,v_2)\longrightarrow(v_2,v_4)\) as a
straight edge-geodesic move.  At the same time, the ordered
sequence
$$
(v_0,v_1,v_2)
\longrightarrow
(v_1,v_2,v_3)
\longrightarrow
(v_2,v_3,v_4)
$$
is the basic local segment of a chamber geodesic.  Thus edge and chamber
geodesics arise from the same local incidence configuration.

The local hypotheses needed to make this continuation consistent are
encoded in the notion of a \emph{weak $A_2$ link}: a bipartite incidence
graph in which induced hexagons serve as apartments and the axioms in
Section~\ref{s:local_geometry} hold.  We study finite colored triangular
complexes $X$ whose vertex links have this form.  Like a
triangulated manifold specified by the topology of its local links, the
present class is specified by local residue conditions. Opposition, successor relations, and
geodesic motion are defined directly from the finite local incidence
data.

This class includes classical $\mathrm{PGL}_3$ building quotients and
quotients of a single apartment. Zeta identities for the latter were
studied by Kang, Li, and Wang in~\cite{KangLiWang2018}.
It also includes compatible strip
gluings and other finite constructions whose zeta data can be defined
without specifying an ambient building, field, lattice, or group action.
The glued examples later in the paper illustrate this flexibility; they
are examples rather than a classification of all complexes satisfying
the local axioms.

The edge and chamber zeta functions of $X$ are defined by the
Euler products
$$
Z_1(X,u)
=
\prod_{[P]\in\mathcal P_1}
\bigl(1-u^{\ell_c(P)}\bigr)^{-1},
\qquad
Z_2(X,u)
=
\prod_{[P]\in\mathcal P_2}
\bigl(1-u^{\ell_c(P)}\bigr)^{-1},
$$
where $\mathcal P_1$ and $\mathcal P_2$ are the sets of cyclic equivalence classes
of primitive closed edge and chamber geodesics, respectively.  Their logarithmic derivatives enumerate closed edge and chamber
geodesics in $X$ of each length, which are the traces of
the powers of the respective weighted successor operators; consequently,
each zeta function is the reciprocal of $\det(\Id-A(u))$ for its weighted
successor operator $A(u)$, so that $Z_1(X,u)$ and $Z_2(X,u)$ are
reciprocal determinants.

The third factor is also determined locally. At each vertex, the
incidence data define a finite augmented local poset.  The M\"obius
function of this poset assembles into a weighted vertex operator
$\widehat A(u)$, whose determinant is the \emph{M\"obius vertex
polynomial} $M(X,u)$ --- a rank-two analogue of the vertex factor
$\det(\Id-Au+(D-\Id)u^2)$ in the classical Ihara--Bass identity for
graphs, where $-1$ and $D-\Id$ are themselves the M\"obius values of
the simpler rank-one local poset at each vertex. When $X$ is a finite
$\mathrm{PGL}_3$ quotient,
$M(X,u)^{-1}$ recovers the usual vertex $L$--factor; see
Remark~\ref{rem:pgl3_mobius_values}. This M\"obius-theoretic viewpoint
on the vertex $L$--factor was first obtained in the setting of general
$\mathrm{PGL}_n$ \cite{KangYu2026}; the present paper develops the same
idea for the purely combinatorial affine $\widetilde A_2$ setting
treated here. Our main theorem is the identity
$$
\frac{Z_1(X,u)}{Z_2(X,-u)}
=
(1-u^3)^{\chi(X)}M(X,u)^{-1},
$$
where \(\chi(X)=|X^{(0)}|-|X^{(1)}|+|X^{(2)}|\) is the ordinary Euler characteristic.
Thus the same three-factor shape that appears in the classical
$\mathrm{PGL}_3$ theory persists in a purely finite combinatorial
setting.

Our proof adapts the cohomological operator framework developed in
\cite{KangLiArtinPGL3} for Artin $L$--functions on finite quotients of the
$\mathrm{PGL}_3$ building.  We construct a pointed complex
\[
C_0(X)\xrightarrow{d_0}C_1(X)\xrightarrow{d_1}C_2(X),
\qquad d_1d_0=0,
\]
together with auxiliary degree-lowering incidence maps $\delta_1$ and
$\delta_2$.  The corresponding degree operators are the components of
\[
\Phi=(1-u^3)\Id+d\delta+\delta d,
\]
so that $\Phi$ is cochain-homotopic to multiplication by $1-u^3$.  The
new point here is that every operator and every local cancellation is
reconstructed from the weak $A_2$ link axioms, without using an ambient
group, Hecke algebra, or representation theory.  For the final determinant
calculation we keep the homotopy interpretation but do not pass to
cohomology: in this three-term complex the alternating cancellation follows
directly from $d_1d_0=0$ and Sylvester's identity
$\det(\Id+UV)=\det(\Id+VU)$.  All operators have entries in $\mathbb Z[u]$,
and the whole argument takes place over the formal power series ring
$\mathbb Z[[u]]$.

This paper is organized as follows.
The local geometry is developed in Section~\ref{s:local_geometry}.
Sections~\ref{s:zeta_functions} and~\ref{s:vertex_mobius_polynomial}
define the successor relations and the three determinant factors.
We prove the main identity in Section~\ref{s:pointed_operators} and
illustrate it with explicit examples in Section~\ref{s:example_verifications}.

\section*{Acknowledgments}
The authors thank Wen-Ching Winnie Li for helpful discussions and valuable comments.

\section{Local Geometry}\label{s:local_geometry}
Throughout this paper, $X$ denotes a pure two-dimensional abstract simplicial complex; it is assumed finite from Section~\ref{s:zeta_functions} on. Write $X^{(k)}$ for the set of $k$-dimensional faces of $X$; elements of $X^{(0)}$, $X^{(1)}$, and $X^{(2)}$ are called \emph{vertices}, 
\emph{edges}, and 
\emph{chambers}, respectively. 

We do not assume that $X$ is the flag completion of its colored 1-skeleton. That is, a colored 3-cycle in $X^{(1)}$ constitutes a chamber if and only if it is explicitly specified as an element of $X^{(2)}$. In particular, when we glue complexes along a common subcomplex, the chambers of the resulting complex are the inherited chambers unless explicitly stated otherwise.

In this section, we introduce the local combinatorial framework for abstract affine complexes $X$ of type $\tilde{A}_2$. The local geometry at each vertex is governed by a bipartite link graph covered by hexagonal apartments, from which all notions of opposition, successor relations, and geodesic motion are constructed.

\subsection{Colored triangular complexes and pointed faces}
A \emph{coloring} of $X$ is a map $\tau \colon X^{(0)} \longrightarrow \mathbb{Z}/3\mathbb{Z}$ such that adjacent vertices have distinct colors; that is, $\{x,y\} \in X^{(1)} \implies \tau(x) \neq \tau(y)$. A simplicial complex equipped with such a coloring is called a \emph{colored complex}, and since $X$ is pure two-dimensional, we refer to it as a \emph{colored triangular complex}.

For vertices $x, y \in X^{(0)}$, the \emph{color difference} from $x$ to $y$ is defined as 
\[
\tau_d(x,y) = \tau(y) - \tau(x) \pmod 3 \in \{0, 1, 2\}.
\]
If $\{x,y\} \in X^{(1)}$, then $\tau_d(x,y) \in \{1,2\}$.

 Combinatorial geodesic dynamics are naturally formulated in terms of pointed faces, which we introduce as follows. For an edge $e = \{x,y\} \in X^{(1)}$, the pair $(e, x)$, also denoted by $(x, y)$, is called a \emph{pointed edge}. The pointed edge $(x,y)$ is of \emph{type $d$} if $\tau_d(x,y) = d$, where $d \in \{1,2\}$. Its reversal is $(y,x)$, which has distinguished vertex $y$ and type $3-d$. 

For a chamber $c = \{x,y,z\} \in X^{(2)}$, the pair $(c, x)$ is called a \emph{pointed chamber} with distinguished vertex $x$. Since the coloring is proper, the vertices of $c$ admit a unique ordering $(x,y,z)$ beginning with $x$ such that 
\[\tau_d(x,y) = \tau_d(y,z) = \tau_d(z,x) = 1.\] 
We also denote $(c,x)$ by $(x,y,z)$.

Note that for $k=1, 2$, each $k$-face gives rise to exactly $k+1$ pointed $k$-faces, corresponding to the $k+1$ choices of distinguished vertex.

\subsection{Weak \texorpdfstring{$A_2$}{A2} links}

We now turn from the two-dimensional complex $X$ to its local structure
at each vertex.

\begin{definition}[Link of a vertex]
Let $x\in X^{(0)}$. The \emph{link} of $x$, denoted by
$\mathrm{Lk}(x)$, is the graph defined as follows:
\begin{enumerate}
\item the vertices of $\mathrm{Lk}(x)$ are the vertices
$y\in X^{(0)}\setminus\{x\}$ such that $\{x,y\}\in X^{(1)}$;
\item two vertices $y,z$ of $\mathrm{Lk}(x)$ are joined by an edge if
and only if $\{x,y,z\}\in X^{(2)}$.
\end{enumerate}
\end{definition}

If $X$ is a colored triangular complex, then $\mathrm{Lk}(x)$ is
bipartite. In the classical setting of affine buildings associated with
$\mathrm{PGL}_3$, such a link is a
spherical building of type $A_2$, whose apartments are hexagons.
Motivated by this local structure, we introduce the following
generalization.

In a bipartite graph, an induced $6$--cycle is called a
\emph{hexagonal apartment}. Throughout this paper, we use the
\emph{complete} apartment system: every induced $6$--cycle is regarded
as an apartment.

\begin{definition}[Weak $A_2$ link]
\label{def:abstract_spherical_A2}
A \emph{weak $A_2$ link} is a connected bipartite graph
$Y$ satisfying the following conditions.
\begin{enumerate}
\item
Every edge of $Y$ is contained in at least one hexagonal apartment.
\item
The graph $Y$ contains no cycles of length shorter than $6$.
Equivalently, $Y$ has girth $6$.
\item
Any two edges of $Y$ are contained in a common hexagonal apartment.
\end{enumerate}
\end{definition}

For a weak $A_2$ link $Y$, we write
$Y=(P\sqcup L,E)$, where $P\sqcup L$ is its vertex bipartition and
$E$ is its edge set.

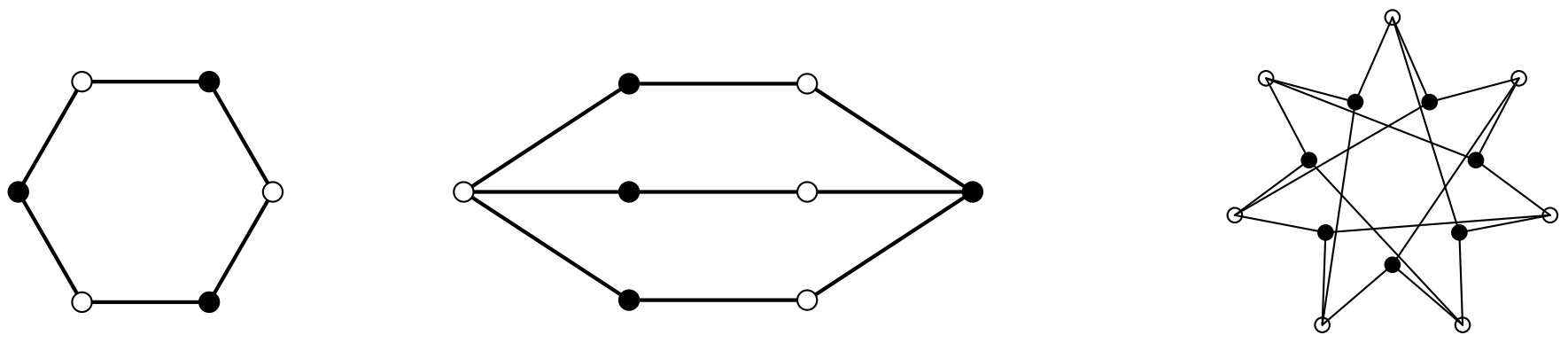
\begin{figure}[ht]
\centering
\begin{tikzpicture}[scale=0.95, line join=round, line cap=round]

  \begin{scope}
    \DrawHexLink
    \node[font=\small] at (0,-1.55)
      {the hexagon $C_6$};
  \end{scope}

  \begin{scope}[shift={(2.5,0)}]
    \DrawThetaLink
    \node[font=\small] at (2,-1.55)
      {the theta graph $\Theta_{3,3,3}$};
  \end{scope}

  \begin{scope}[shift={(9.8,.1)}, scale=0.82]

    \foreach \i in {0,...,6} {
      \coordinate (p\i) at ({90-360*\i/7}:1.55);
    }

    \foreach \i in {0,...,6} {
      \coordinate (l\i) at ({90-360*\i/7+180/7}:0.82);
    }

    \foreach \i in {0,...,6} {
      \filldraw[fill=white] (p\i) circle (2pt);
      \filldraw[fill=black] (l\i) circle (2pt);
    }

    \foreach \i in {0,...,6} {
      \pgfmathtruncatemacro{\j}{mod(\i+1,7)}
      \pgfmathtruncatemacro{\k}{mod(\i+3,7)}
      \draw (p\i) -- (l\i);
      \draw (p\i) -- (l\j);
      \draw (p\i) -- (l\k);
    }

    \node[font=\small] at (0,-2.05)
      {the incidence graph of $\mathbb{P}^2(\mathbb{F}_2)$};

  \end{scope}

\end{tikzpicture}

\caption{Three weak $A_2$ links}
\label{fig:weak-A2-examples}
\end{figure}

The definition imposes no thickness condition.  Thus the thin hexagon
$C_6$ and the graph $\Theta_{3,3,3}$ pictured in Figure~\ref{fig:weak-A2-examples} below are both allowed. 

\begin{example}[The hexagon]
The graph $C_6$ is the thin local model.  Its only hexagonal apartment is
the whole graph, so every edge lies in a hexagonal apartment and any two
local edges lie in a common hexagonal apartment.
\end{example}

\begin{example}[The theta graph $\Theta_{3,3,3}$]
The graph $\Theta_{3,3,3}$ consists of three length-three paths with
common endpoints.  Any two of the three paths form a hexagonal apartment.
Hence every edge lies in a hexagonal apartment, and any two local edges
lie in a common hexagonal apartment.  This is the smallest non-hexagonal weak $A_2$ link.  More generally, for any $r\geq 2$, the graph
$\Theta_{3,\ldots,3}$ consisting of $r$ length-three paths with common
endpoints is a weak $A_2$ link; the case $r=2$ is the hexagon $C_6$.
\end{example}

\begin{example}[The projective-plane incidence graph]
\label{ex:projective-plane-link}
Let $Y_q$ be the incidence graph of the projective plane
$\mathbb{P}^2(\mathbb{F}_q)$.  Thus the two parts of the bipartition are
the points and lines of $\mathbb{P}^2(\mathbb{F}_q)$, with a point
adjacent to a line precisely when the point lies on the line.
This is the incidence graph of the spherical building of type $A_2$
over $\mathbb{F}_q$, and it is the motivating example for
Definition~\ref{def:abstract_spherical_A2}.  Its apartments are the
hexagons arising from triples of points and lines in general position.
In particular, every edge is contained in a hexagonal apartment.

Moreover, the incidence graph has girth $6$: two distinct points lie
on a unique common line, and dually two distinct lines meet in a unique
point, so there are no $4$-cycles.  Finally, the building apartment
axiom implies that any two edges are contained in a common apartment.
Hence $Y_q$ is a weak $A_2$ link. The case $q=2$ is shown in Figure~\ref{fig:weak-A2-examples}.
\end{example}

\begin{definition}
Two vertices $y,z$ in a weak $A_2$ link $Y$ are said to be \emph{opposite}
if they have graph distance $3$. Equivalently, in any hexagonal apartment
containing $y$ and $z$, they are furthest apart.
\end{definition}

The following lemma records some properties of opposite vertices. 

\begin{lemma}[Adjacency and opposition in a weak $A_2$ link]
\label{lem:spherical_A2_basic}
Let $Y=(P\sqcup L,E)$ be a weak $A_2$ link.
Then:
\begin{enumerate}
\item any two distinct vertices in $P$ (resp.\ $L$) have a unique
common neighbor in $L$ (resp.\ $P$);
\item given $p\in P$ and $\ell\in L$, the vertices $p$ and $\ell$ are
opposite if and only if they are not adjacent.
\end{enumerate}
In particular, the notion of opposition is independent of the hexagonal
apartment chosen to witness it.
\end{lemma}

\begin{proof}
We prove the statements for vertices in $P$; the argument for $L$ is similar.  Let
$p,p'\in P$ be distinct.  Since $Y$ is connected,  choose an edge through $p$
and an edge through $p'$.  By the apartment axiom, these two edges are contained in
a hexagonal apartment.  In a hexagonal apartment, two distinct
vertices in $P$ 
share a common neighbor in $L$; this proves
existence.  If two different vertices in $L$ were adjacent with both $p$
and $p'$, they would form a $4$--cycle, contradicting the girth-$6$
assumption.  This proves uniqueness, and hence the first statement.

Now let $p\in P$ and $\ell\in L$.  If $p$ and $\ell$ are adjacent, then
their graph distance is $1$, so they are not opposite by definition.  Conversely, suppose that $p$ and $\ell$ are not adjacent.
Choose an edge containing $p$ and an edge containing $\ell$.
By the apartment axiom, these two edges lie in a common hexagonal apartment.
Since the apartment is bipartite, two vertices of different types have
distance either $1$ or $3$ within the apartment.  As $p$ and $\ell$ are
not adjacent, their distance is $3$, and hence they are opposite.
This proves the second statement.
\end{proof}

\subsection{Abstract affine complexes of type
\texorpdfstring{$\widetilde A_2$}{affine A2}}

We now impose the weak $A_2$ link condition at every vertex of a colored
triangular complex.  The coloring separates each link graph into local
points and local lines.

\begin{definition}[Abstract affine complex of type $\widetilde A_2$]
\label{def:abstract_affine_A2}
An \emph{abstract affine complex of type $\widetilde A_2$} is a pure
colored triangular complex $(X,\tau)$ satisfying the following local
condition:  For every vertex $x\in X^{(0)}$, let
\(P_x=\{y\in \mathrm{Lk}(x)^{(0)}\mid \tau_d(x,y)=1\}\) and
\(L_x=\{y\in \mathrm{Lk}(x)^{(0)}\mid \tau_d(x,y)=2\}\).  Then the link
$\mathrm{Lk}(x)$, regarded as the bipartite graph
\(\mathrm{Lk}(x)=(P_x\sqcup L_x,E_x)\), is a weak $A_2$ link.
\end{definition}

 The above definition is entirely local. All incidence structures and local geometric data required in our constructions take place within the vertex links $\mathrm{Lk}(x)$; no global affine building structure, global system of apartments, or ambient group action is assumed. This local perspective is the two-dimensional analogue of studying arbitrary finite graphs without requiring them to be quotients of rank-one buildings.  While passing to the universal cover $\widetilde{X}$ is a standard and useful tool---much as in graph theory or topological manifold theory---we do not require $\widetilde{X}$ to be an affine building of type $\widetilde{A}_2$.

\subsection{Examples of abstract affine complexes of type
\texorpdfstring{$\widetilde A_2$}{affine A2}}
\subsubsection{The standard affine apartment}\label{s:standard_affine_apartment}

Before giving examples of finite complexes, we fix the apartment model used in the
figures.  Let \(\mathcal A\) be the triangular tiling with vertex set
\(\mathbb Z^2\).  We join \((i,j)\) to its neighbors in the directions
\((1,0)\), \((0,1)\), and \((1,1)\), and their negatives.  We color the vertices by
\(\operatorname{col}(i,j)=i+j\pmod 3\).
With this convention every chamber has one vertex of each color.
Translations \((i,j)\mapsto(i+a,j+b)\) preserve the coloring exactly
when \(a+b\equiv 0\pmod 3\).  Quotients of this apartment by such
translations give flat triangular tori; adding a glide reflection gives
Klein bottle quotients.  The gluing examples below are obtained by
identifying compatible annular strips in such apartment quotients.

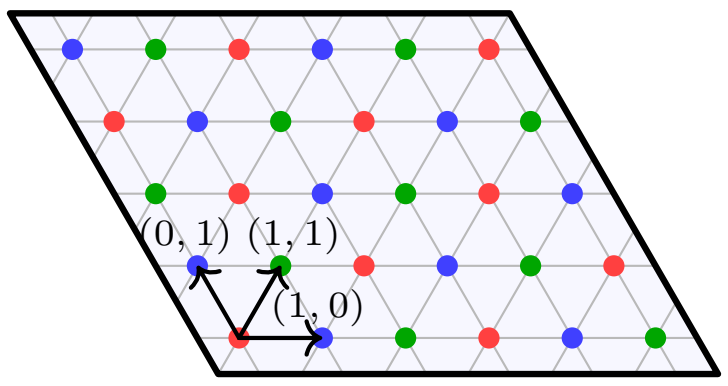
\begin{figure}[ht]
\centering
\begin{tikzpicture}[scale=0.95, line join=round, line cap=round]
\DrawAffineApartmentPatch
\end{tikzpicture}
\caption{A finite patch of the standard affine apartment $\mathcal A$ with the three edge directions and the color pattern.}
\label{fig:standard-affine-apartment}
\end{figure}

\subsubsection{Apartment quotients}

The next examples are finite quotients of the standard affine apartment.
In both cases the acting group acts freely and preserves types, so every
vertex link in the quotient is the hexagon $C_6$.  We therefore specify
only the quotient data and the features needed later.  Compatible
annular strips in these quotients will subsequently provide the building
blocks for the gluing examples.

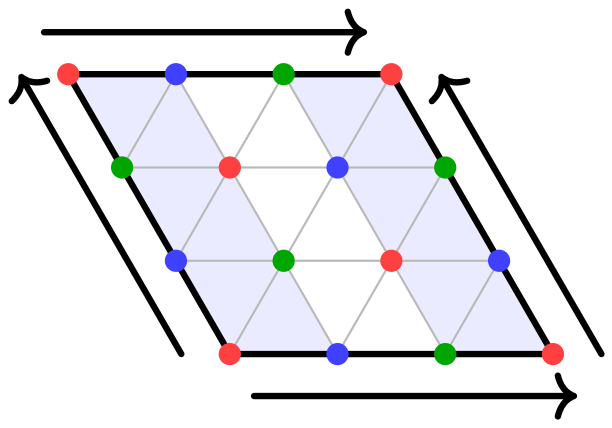
\begin{figure}[ht]
\centering
\begin{tikzpicture}[scale=0.95, line join=round, line cap=round]

\DrawCorrectTorusCopy{(0,0)}{blue}{T}{white!0}

\end{tikzpicture}
\caption{The minimal apartment torus.}
\label{fig:minimal-apartment-torus}
\end{figure}

\begin{example}[Minimal apartment torus]\label{ex:minimal-apartment-torus}
The minimal apartment torus is the quotient of $\mathcal A$ by the
translation lattice generated by $(3,0)$ and $(0,3)$.  Its $3\times3$
fundamental parallelogram and side identifications are shown in
Figure~\ref{fig:minimal-apartment-torus}.  Every vertex link is $C_6$.

Its numerical data are \(|X^{(0)}|=9\), \(|X^{(1)}|=27\),
\(|X^{(2)}|=18\), and \(\chi(X)=0\).
\end{example}

\begin{figure}[ht]
\centering
\begin{tikzpicture}[scale=0.95, line join=round, line cap=round]

\DrawKleinBottleCopy{(0,0)}{blue}{gray!22}

\end{tikzpicture}
\caption{The torus double cover of the Klein bottle quotient.}
\label{fig:klein-bottle-quotient}
\end{figure}
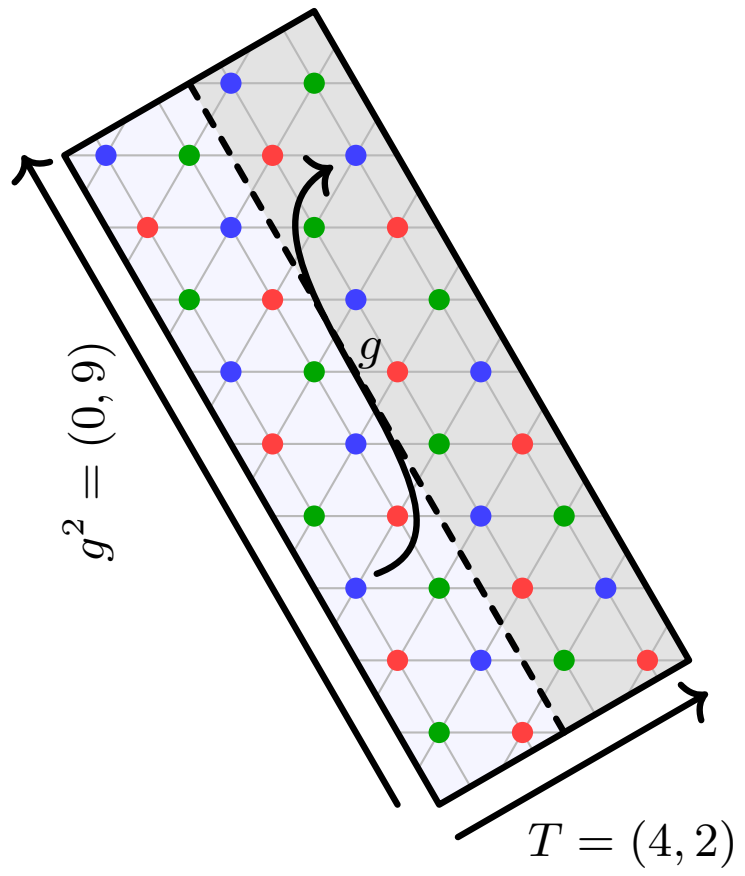

\begin{example}[Klein bottle quotient]\label{ex:klein-bottle-quotient}
Let \(T(i,j)=(i+4,j+2)\) and
\(g(i,j)=(-i-1,\,-i+j+4)\).  Both maps preserve the colored triangular
tiling.  Moreover, \(g^2(i,j)=(i,j+9)\) and \(gTg^{-1}=T^{-1}\).
The translation subgroup $\langle T,g^2\rangle$ gives the torus double
cover
\[
\mathbb T_{\mathrm{Kl}}
=
\langle T,g^2\rangle\backslash\mathcal A
\]
with periods $(4,2)$ and $(0,9)$.  The glide reflection $g$ descends to
a fixed-point-free involution $\bar g$, and the Klein bottle quotient is
\[
X_{\mathrm{Kl}}
=
\mathbb T_{\mathrm{Kl}}/\langle\bar g\rangle
=
\langle T,g\rangle\backslash\mathcal A.
\]
Figure~\ref{fig:klein-bottle-quotient} shows the torus double cover.
Every vertex link is $C_6$, and \(|X^{(0)}|=18\),
\(|X^{(1)}|=54\), \(|X^{(2)}|=36\), and \(\chi(X)=0\).
The glide reflection will make this example useful for detecting the
sign in $Z_2(X,-u)$.
\end{example}

\subsubsection{Gluing along strips}\label{s:strip_gluing}
We now describe a local gluing operation that builds new abstract affine
complexes of type $\widetilde A_2$ from old ones, by identifying an
annular strip in each; the resulting constructions supply the glued
examples used later in the paper.

\begin{definition}[Roots in weak $A_2$ links]
Let $Y$ be a weak $A_2$ link. A \emph{root} in $Y$ is a path
of length $3$.
A root $R$ is \emph{admissible} if, for every edge $e$ of $Y$,
there is a hexagonal apartment of $Y$ containing both $R$ and $e$.
\end{definition}

\begin{lemma}[Root gluing of weak $A_2$ links]
\label{lem:root_gluing}
Let $Y_1$ and $Y_2$ be weak $A_2$ links, and let
$R_i\subset Y_i$, $i=1,2$, be admissible roots.
If $\psi:R_1\longrightarrow R_2$ is a graph isomorphism,
then the graph $Y=Y_1\cup_{\psi}Y_2$ obtained by identifying
the two roots via $\psi$ is again a weak $A_2$ link.
\end{lemma}

\begin{proof}
Let $R$ denote the common image of $R_1$ and $R_2$ in $Y$.
The graph $Y$ is connected because the connected graphs $Y_1$
and $Y_2$ are glued along the nonempty path $R$.
After interchanging the two parts of the bipartition of $Y_2$
if necessary, $\psi$ preserves the bipartitions, so $Y$ is
also bipartite.

We first show that any two edges $e$ and $f$ of $Y$ lie in
a common hexagonal apartment. If both lie in the same $Y_i$,
this follows from the weak $A_2$ property of $Y_i$.
Otherwise, we may assume that
$e\in E(Y_1)\setminus E(R_1)$ and
$f\in E(Y_2)\setminus E(R_2)$.
By admissibility, choose hexagonal apartments
$A_1\subset Y_1$ containing $R_1$ and $e$, and
$A_2\subset Y_2$ containing $R_2$ and $f$.
In $Y$, the complementary paths of length $3$ of $R$ in $A_1$ and in
$A_2$ meet only at the endpoints of $R$, so their union is a hexagon
containing both $e$ and $f$.

It remains to verify the girth condition.
Any cycle contained in either $Y_i$ has length at least $6$.
Otherwise, a cycle $C$ contains, for each $i=1,2$, a path
$P_i$ in $Y_i$ with edges outside $R$, distinct endpoints on $R$,
and interior disjoint from $R$.
Let $\ell_i$ be its length and $d_i$ the distance between its
endpoints along $R$. Since $P_i$ and the corresponding subpath
of $R$ form a cycle in $Y_i$,
$$
\ell_i+d_i\geq 6,
\qquad 1\leq d_i\leq 3.
$$
Thus each $P_i$ has length at least $3$. These paths are
edge-disjoint, so $C$ has length at least $6$.
Therefore $Y$ has no cycles of length less than $6$.
In particular, every hexagon constructed above is induced.
Hence $Y$ is a weak $A_2$ link.
\end{proof}

\begin{definition}[Admissible strip]
Let $X$ be an abstract affine complex of type $\widetilde A_2$.  A
colored subcomplex $S\subset X$ is called an \emph{admissible strip} if
for every vertex $x\in S$, the link \(\operatorname{Lk}_S(x)\) is an
admissible root in the weak $A_2$ link \(\operatorname{Lk}_X(x)\).
If, furthermore, $S$ is homeomorphic to an annulus, we call it a \emph{closed
admissible strip}.
\end{definition}

\begin{proposition}[Gluing principle for abstract affine complexes]
\label{prop:strip_gluing}
Let $X_1$ and $X_2$ be abstract affine complexes of type
$\widetilde A_2$, and let \(S_i\subset X_i, i=1,2,\) be admissible strips.
Suppose there is a color-preserving simplicial
isomorphism \(\varphi:S_1\longrightarrow S_2\). Assume that, after
the vertex identification, the images of $X_1$ and $X_2$ intersect
exactly in the common strip (as holds if either $S_i$ is a full
subcomplex of $X_i$). Let
\(X=X_1\cup_{\varphi}X_2\) be the complex obtained by
identifying $S_1$ with $S_2$.  Then $X$ is again
an abstract affine complex of type $\widetilde A_2$.
\end{proposition}

\begin{proof}
The color-preserving simplicial gluing preserves the proper coloring and
purity of the triangular chamber complex, since the chambers are exactly
the inherited chambers.  It remains to verify the local link condition.
If $x\in X$ is not on the common image of the
glued strip, then a neighborhood of $x$ is unchanged, so
\(\operatorname{Lk}_X(x)\) is the same weak $A_2$ link as before.

Now suppose $x$ lies on the common image of the glued strip, and let
$x_i\in S_i$ be its preimages.  By the intersection hypothesis, the link of $x$ in the glued complex is
the pushout
\[
\operatorname{Lk}_X(x)
=
\operatorname{Lk}_{X_1}(x_1)
\cup_{\operatorname{Lk}_{S_1}(x_1)\cong
      \operatorname{Lk}_{S_2}(x_2)}
\operatorname{Lk}_{X_2}(x_2).
\]
By definition of admissible strip, the common subgraph
\(\operatorname{Lk}_{S_i}(x_i)\) is an admissible root in each ambient
link.  Therefore Lemma~\ref{lem:root_gluing} implies that
\(\operatorname{Lk}_X(x)\) is again a weak $A_2$ link.  Hence every
vertex of $X$ has weak $A_2$
link, so $X$ is an abstract affine complex of type $\widetilde A_2$.
\end{proof}

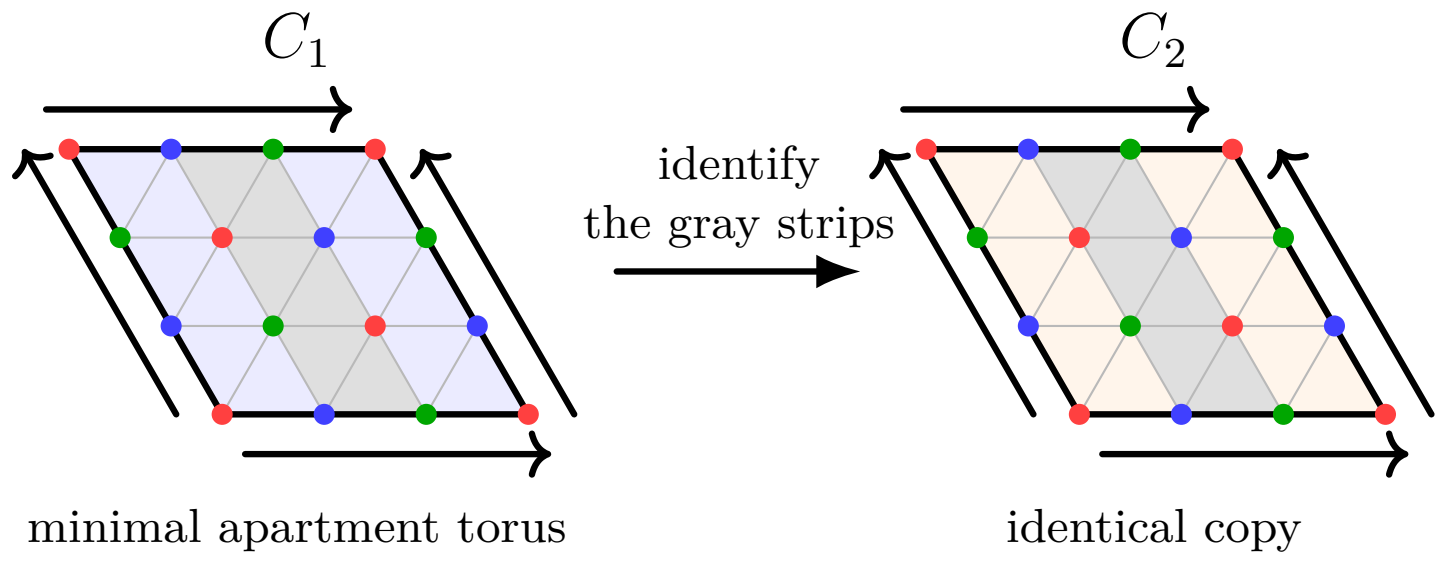
\begin{figure}[ht]
\centering
\begin{tikzpicture}[scale=0.90, line join=round, line cap=round]

  \DrawCorrectTorusCopy{(0,0)}{blue}{A}{gray!25}
  \node[font=\large] at (0.55,2.75) {$C_1$};
  \node[font=\small] at (0.55,-0.85) {minimal apartment torus};

  \DrawCorrectTorusCopy{(6.3,0)}{orange}{B}{gray!25}
  \node[font=\large] at (6.85,2.75) {$C_2$};
  \node[font=\small] at (6.85,-0.85) {identical copy};

  \draw[-{Latex[length=3mm]}, very thick] (2.90,1.05) -- (4.70,1.05);
  \node[align=center, font=\small] at (3.80,1.60)
    {identify\\the gray strips};

\end{tikzpicture}
\caption{Two minimal apartment tori glued along a common annular strip.}
\label{fig:two-torus-gluing}
\end{figure}

\begin{example}[Two tori glued along a strip]\label{ex:two-torus-gluing}
Take two copies $C_1$ and $C_2$ of the minimal apartment torus and
identify the gray annular strips shown in
Figure~\ref{fig:two-torus-gluing} as full colored subcomplexes, using the
inherited chamber set. Proposition~\ref{prop:strip_gluing} applies.
Vertices off the common strip have link $C_6$,
and vertices on it have link $\Theta_{3,3,3}$.  The resulting abstract
affine complex has \(|X^{(0)}|=12\), \(|X^{(1)}|=42\),
\(|X^{(2)}|=30\), and \(\chi(X)=0\).
\end{example}

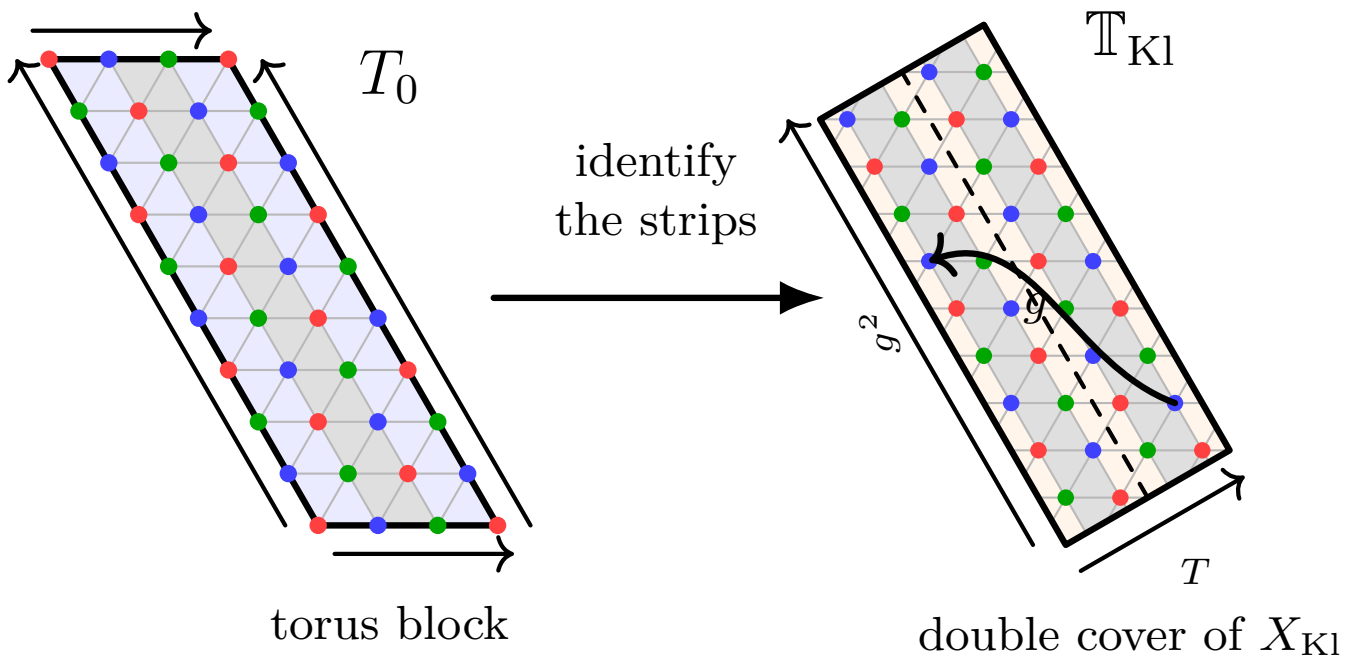
\begin{figure}[ht]
\centering
\begin{tikzpicture}[scale=0.86, line join=round, line cap=round]

  \DrawTorusNineStripBlock{(2,0)}{blue}{gray!25}
  \node[font=\large] at (2.55,3.45) {$T_0$};
  \node[font=\small] at (2.55,-0.75) {torus block};

  \DrawKleinDoubleCoverWithLiftedStrips{(8.8,-0.15)}{orange}{gray!25}
  \node[font=\large] at (8.25,3.75) {$\mathbb T_{\mathrm{Kl}}$};
  \node[font=\small] at (8.25,-0.85) {double cover of $X_{\mathrm{Kl}}$};

  \draw[-{Latex[length=3mm]}, very thick] (3.35,1.75) -- (5.90,1.75);
  \node[align=center, font=\small] at (4.60,2.55)
    {identify\\the strips};

\end{tikzpicture}
\caption{A torus and Klein bottle quotient glued along compatible
annular strips.}
\label{fig:torus-klein-gluing}
\end{figure}

\begin{example}[A torus and a Klein bottle glued along a strip]
\label{ex:torus-klein-gluing}
Let \(T_0=\langle (3,0),(0,9)\rangle\backslash\mathcal A\) be a torus
quotient of the standard apartment.  It contains a vertical
one-chamber-wide annular strip with $18$ chambers.  In the Klein bottle
$X_{\mathrm{Kl}}$ from Example~\ref{ex:klein-bottle-quotient}, the strip
centered on the glide axis projects to a M\"obius band with $9$ chambers
and is not used here.  A vertical strip away from that axis is paired
with a disjoint strip in the torus double cover, and the pair projects to
a two-sided annular strip with $18$ chambers.

The strip in $T_0$ is a full subcomplex, so the intersection hypothesis
of Proposition~\ref{prop:strip_gluing} holds. We identify the two annular strips
color-preservingly. On a chosen lift,
the identification is
\[
(1,j)\longmapsto(0,j+1),\qquad
(2,j)\longmapsto(1,j+1),
\qquad j\in\mathbb Z/9\mathbb Z;
\]
the shift by $1$ preserves the coloring.  With the inherited chamber
set, let \(X=T_0\cup_S X_{\mathrm{Kl}}\).
Figure~\ref{fig:torus-klein-gluing} shows the construction.  The torus
block is drawn directly.  The Klein bottle block is drawn through its
torus double cover: the two gray lifted strips are exchanged by the glide
involution and project to the annular strip used for the gluing.  There
are $9$ vertices with link $C_6$ and $18$ vertices on the common strip
with link $\Theta_{3,3,3}$.  The numerical data are
\(|X^{(0)}|=27\), \(|X^{(1)}|=99\), \(|X^{(2)}|=72\), and
\(\chi(X)=0\).
\end{example}

\section{Edge and Chamber Zeta Functions}\label{s:zeta_functions}

The main purpose of this section is to introduce edge and chamber geodesics in $X$ and systematically study their associated zeta functions. We begin by reviewing local opposition and establishing the successor data for pointed edges and chambers. Building on these local structures, we define the edge and chamber zeta functions $Z_1(X, u)$ and $Z_2(X, u)$.  Throughout this section, $X$ is finite.

We describe these zeta functions from several equivalent viewpoints: as Euler products over primitive closed geodesics, as exponential series via based combinatorics, and as reciprocal determinants $\det(\Id - A(u))^{-1}$ of weighted adjacency operators on pointed cells.

\subsection{Successor relation}\label{s:opposite_successor}
Let $X$ be an abstract affine complex of type $\widetilde A_2$.
Since every link $\mathrm{Lk}(x)$ is a weak $A_2$ link,
opposition is canonically defined in each link.
We now define the successor relation on pointed faces;
it plays the role of non-backtracking in ordinary graph theory.

\begin{definition}[Successor]\label{def:successor_relation}
A \emph{successor} of a pointed edge $(v_0,v_1)$ is a pointed edge
$(v_1,v_2)$ such that the edges $\{v_0,v_1\}$ and $\{v_1,v_2\}$
are opposite at $v_1$, that is, the vertices $v_0$ and $v_2$
of $\mathrm{Lk}(v_1)$ are opposite.

A \emph{successor} of a pointed chamber $(v_0,v_1,v_2)$ is a pointed
chamber $(v_1,v_2,v_3)$ such that the chambers
$\{v_0,v_1,v_2\}$ and $\{v_1,v_2,v_3\}$ are distinct and share the
common edge $\{v_1,v_2\}$.
\end{definition}

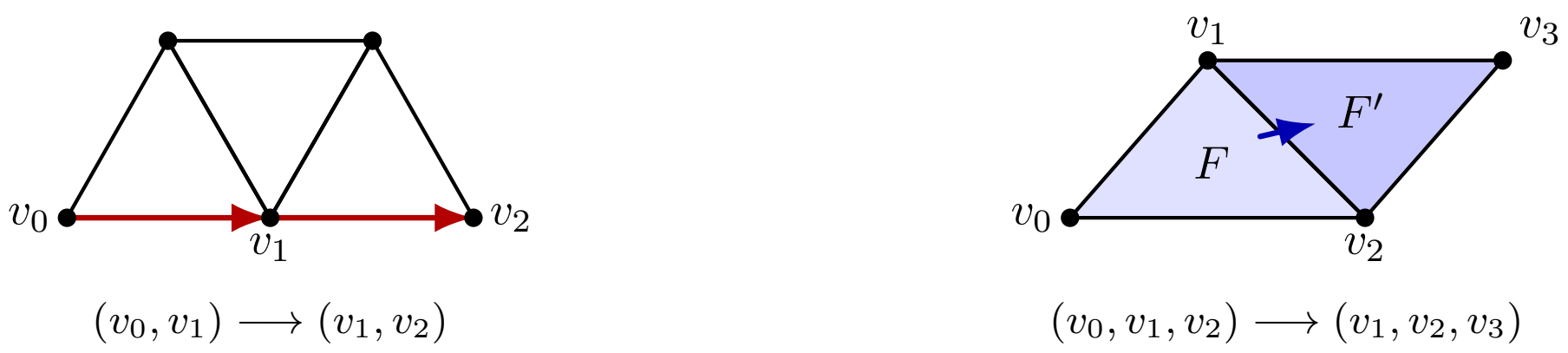
\begin{figure}[ht]
\centering
\begin{tikzpicture}[
    scale=1.0,
    line join=round,
    line cap=round,
    >=Latex
]

  \begin{scope}[shift={(-3.7,0)}]

    \coordinate (v1) at (0,0);
    \coordinate (v0) at (-1.55,0);
    \coordinate (u1) at (-0.78,1.35);
    \coordinate (u2) at (0.78,1.35);
    \coordinate (v2) at (1.55,0);

    \draw[thick] (v1) -- (v0) -- (u1) -- cycle;
    \draw[thick] (v1) -- (u1) -- (u2) -- cycle;
    \draw[thick] (v1) -- (u2) -- (v2) -- cycle;

    \draw[very thick,red!70!black,->] (v0) -- (v1);
    \draw[very thick,red!70!black,->] (v1) -- (v2);

    \foreach \p in {v0,v1,v2,u1,u2} {
      \fill (\p) circle (2pt);
    }

    \node[left] at (v0) {$v_0$};
    \node[below] at (v1) {$v_1$};
    \node[right] at (v2) {$v_2$};

    \node[font=\small] at (0,-0.8)
      {$(v_0,v_1)\longrightarrow (v_1,v_2)$};

  \end{scope}

  \begin{scope}[shift={(3.9,0)}]

    \coordinate (v0) at (-1.5,0);
    \coordinate (v1) at (-0.45,1.2);
    \coordinate (v2) at (0.75,0);
    \coordinate (v3) at (1.8,1.2);

    \fill[blue!12] (v0) -- (v1) -- (v2) -- cycle;
    \fill[blue!22] (v1) -- (v2) -- (v3) -- cycle;

    \draw[thick] (v0) -- (v1) -- (v2) -- cycle;
    \draw[thick] (v1) -- (v2) -- (v3) -- cycle;

    \foreach \p in {v0,v1,v2,v3} {
      \fill (\p) circle (2pt);
    }

    \node[left] at (v0) {$v_0$};
    \node[above] at (v1) {$v_1$};
    \node[below] at (v2) {$v_2$};
    \node[above right] at (v3) {$v_3$};

    \node at (-0.42,0.42) {$F$};
    \node at (0.72,0.82) {$F'$};

    \draw[blue!70!black,very thick,->]
      (-0.05,0.62) -- (0.38,0.72);

    \node[font=\small] at (0.15,-0.8)
      {$(v_0,v_1,v_2)\longrightarrow (v_1,v_2,v_3)$};

  \end{scope}

\end{tikzpicture}
\caption{The successor relations for pointed edges and pointed chambers.}
\label{fig:successor-relations}
\end{figure}

\subsection{Edge and chamber geodesics}

Recall that a pointed edge is an ordered pair \((v_0,v_1)\) with
\(\{v_0,v_1\}\in X^{(1)}\), so that
\(\tau_d(v_0,v_1)\in\{1,2\}\). The successor relation on pointed edges
therefore gives a natural notion of directed geodesic motion.

\begin{definition}[Edge geodesic]\label{def: Edge geodesic}
An \emph{edge geodesic} of length $m$ is a sequence of pointed edges in $X$
\[
(v_0,v_1),\,
(v_1,v_2),\,
\dots,\,
(v_{m-1},v_m)
\]
such that, for $1\le i\le m-1$, each pointed edge $(v_i,v_{i+1})$ is a
successor of $(v_{i-1},v_i)$.
\end{definition}

There is an analogous construction for chambers. Recall
that a pointed chamber is a cyclically ordered triple
$(v_0,v_1,v_2)$ satisfying
\[
\tau_d(v_0,v_1)=\tau_d(v_1,v_2)=\tau_d(v_2,v_0)=1.
\]

\begin{definition}[Chamber geodesic]\label{def: Chamber geodesic}
A \emph{chamber geodesic} of length $m$ is a sequence of pointed chambers in $X$
\[
(v_0,v_1,v_2),\,
(v_1,v_2,v_3),\,
\dots,\,
(v_{m-1},v_m,v_{m+1})
\]
such that, for $1\le i\le m-1$, each pointed chamber
$(v_i,v_{i+1},v_{i+2})$ is a successor of
$(v_{i-1},v_i,v_{i+1})$.
\end{definition}

The notions relevant to periodic geodesics are the same in the edge
and chamber settings. A geodesic is \emph{closed} if its first pointed
edge or chamber is a successor of its last one.

Two closed geodesics of the same kind are \emph{cyclically equivalent}
if one is obtained from the other by a cyclic shift. We write $[P]$
for the corresponding cyclic equivalence class.

A closed geodesic $P$ is \emph{primitive} if it is not a nontrivial
power of a shorter closed geodesic. Equivalently, $P$ is \emph{primitive} if there do not exist a
closed geodesic $Q$ and an integer $k\geq 2$ such that
\[
P=Q^k.
\]
Every closed geodesic can be written uniquely as a power of a primitive
one.

The notion of length in Definitions~\ref{def: Edge geodesic} and~\ref{def: Chamber geodesic} counts successor steps. For the construction of zeta functions, we will need the notion of color length introduced below. 
\begin{definition}[Edge color length]
Let
\[
P=
(v_0,v_1),\,
(v_1,v_2),\,
\dots,\,
(v_{m-1},v_m)
\]
be an edge geodesic. Its \emph{color length} is
\[
\ell_c(P)
=
\sum_{i=0}^{m-1}\tau_d(v_i,v_{i+1}).
\]
Along an edge geodesic, all pointed edges have the same type. Hence
either
\[
\tau_d(v_i,v_{i+1})=1
\qquad\text{for all }i,
\]
or
\[
\tau_d(v_i,v_{i+1})=2
\qquad\text{for all }i.
\]
Accordingly,
\[
\ell_c(P)=m
\qquad\text{or}\qquad
\ell_c(P)=2m.
\]
\end{definition}

\begin{definition}[Chamber color length]
Let
\[
P=
(v_0,v_1,v_2),\,
(v_1,v_2,v_3),\,
\dots,\,
(v_{m-1},v_m,v_{m+1})
\]
be a chamber geodesic. Its \emph{color length} is defined analogously by
\[
\ell_c(P)
=
\sum_{i=0}^{m-1}\tau_d(v_i,v_{i+1}).
\]
Along a chamber geodesic, the pointed chambers are ordered so that
\[
\tau_d(v_i,v_{i+1})=1
\qquad\text{for all }i.
\]
Thus every term in the sum is equal to \(1\), and hence
\[
\ell_c(P)=m.
\]
\end{definition}

\subsection{Edge and chamber zeta functions}

Let $\mathcal P_1$ and $\mathcal P_2$ denote the sets of
cyclic equivalence classes of primitive closed edge geodesics and
primitive closed chamber geodesics respectively.

\begin{definition}[Edge zeta function]
The \emph{edge zeta function} of $X$ is the formal power series in $\mathbb Z[[u]]$ given by 
\[
Z_1(X,u)
=
\prod_{[Q]\in\mathcal P_1}
\bigl(1-u^{\ell_c(Q)}\bigr)^{-1}.
\]
\end{definition}

\begin{definition}[Chamber zeta function]
The \emph{chamber zeta function} of $X$ is the formal power series in $\mathbb Z[[u]]$ given by 
\[
Z_2(X,u)
=
\prod_{[Q]\in\mathcal P_2}
\bigl(1-u^{\ell_c(Q)}\bigr)^{-1}.
\]
\end{definition}

Taking the logarithm of $Z_i(X, u)$ in $\mathbb Q[[u]]$ and applying the formula
\[
-\log(1-x)=\sum_{k\geq 1}\frac{x^k}{k},
\]
one obtains
\[
\log Z_i(X,u)
=
\sum_{[Q]\in\mathcal P_i}\sum_{k\geq 1}
\frac{u^{k\ell_c(Q)}}{k}.
\]

Let $\mathcal C_1$ and $\mathcal C_2$ denote the sets of cyclic equivalence classes of all closed edge or chamber geodesics, respectively, including
nonprimitive ones. Recall that for every $[P]\in\mathcal C_i$, the closed geodesic $P$ can be written uniquely as
$P=Q^{\nu(P)}$, where $Q$ is primitive and $\nu(P)\geq 1$.
Since $\ell_c(Q^k)=k\ell_c(Q)$, the pairs $([Q],k)$ with $[Q]\in\mathcal P_i$ and $k\ge 1$ correspond precisely
to the classes $[P]\in\mathcal C_i$, with $\nu(P)=k$. Hence
\begin{equation}\label{eq:all_closed_geodesic_expansion}
\log Z_i(X,u)
=
\sum_{[P]\in\mathcal C_i}
\frac{u^{\ell_c(P)}}{\nu(P)},
\qquad
Z_i(X,u)
=
\exp\!\left(
\sum_{[P]\in\mathcal C_i}
\frac{u^{\ell_c(P)}}{\nu(P)}
\right).
\end{equation}
Thus the zeta functions defined through Euler products over primitive closed geodesics encode information of all closed
geodesics, with a $k$-fold repetition weighted by $1/k$.

There is also a based version of this counting interpretation.  For $i=1, 2$, let
$\mathcal B_{i,r}$ be the set of closed edge or chamber geodesics with a
distinguished initial pointed face and exactly $r$ successor steps, and
put
\[
N_{i,r}(u)
=
\sum_{P\in\mathcal B_{i,r}}u^{\ell_c(P)}.
\]
A primitive cycle with $s$ successor steps has $s$ choices of its
distinguished initial pointed face.  Regrouping its repetitions therefore
gives
\begin{equation}\label{eq:based_closed_geodesic_expansion}
\log Z_i(X,u)
=
\sum_{r\geq 1}\frac{N_{i,r}(u)}{r}.
\end{equation}
For chambers, every successor step has color increment $1$. Let $N_{2,r}$ be the number of elements in $\mathcal B_{2,r}$, then
\[
Z_2(X,u)
=
\exp\!\left(\sum_{r\geq 1}\frac{N_{2,r}}{r}u^r\right).
\]
For edges, the successor relation preserves edge type.  Let
$N_{1,r}^{(a)}$ be the number of elements in $\mathcal B_{1,r}$ of type
$a\in\{1,2\}$, then
\[
Z_1(X,u)
=
\exp\!\left(
\sum_{r\geq 1}
\frac{N_{1,r}^{(1)}u^r+N_{1,r}^{(2)}u^{2r}}{r}
\right).
\]

\subsection{The Hashimoto determinant formulas}

 
Let \(C_1(X)\) and \(C_2(X)\) denote the free $\mathbb Z[[u]]$-modules of
$\mathbb Z[[u]]$-valued functions on the pointed
edges and the pointed chambers of $X$, respectively.
\begin{definition}[Weighted edge adjacency operator]
\label{def:weighted_edge_adjacency}
The \emph{weighted edge adjacency operator}
\(A_E(u):C_1(X)\longrightarrow C_1(X)\) is defined by
\[
(A_E(u)f)(v_0,v_1)
=
\sum_{\substack{(v_1,v_2)\\
(v_1,v_2)\succ(v_0,v_1)}}
u^{\tau_d(v_1,v_2)}
f(v_1,v_2),
\]
where the summation runs over all successors of the pointed edge
$(v_0,v_1)$.
\end{definition}

\begin{definition}[Weighted chamber adjacency operator]
\label{def:weighted_chamber_adjacency}
The \emph{weighted chamber adjacency operator}
\(A_C(u):C_2(X)\longrightarrow C_2(X)\) is defined by
\[
(A_C(u)f)(v_0,v_1,v_2)
=
\sum_{\substack{(v_1,v_2,v_3)\\
(v_1,v_2,v_3)\succ(v_0,v_1,v_2)}}
u^{\tau_d(v_2,v_3)}
f(v_1,v_2,v_3),
\]
where the summation runs over all successors of the pointed chamber
$(v_0,v_1,v_2)$.
\end{definition}

The edge and chamber zeta functions can be viewed as zeta functions of finite weighted
directed graphs constructed in the following way.  The vertices of the graph are the
pointed edges or pointed chambers, and the directed edges are defined by the associated successor transitions.  The weight of a transition is determined
by the color difference of its terminal pointed face.

By construction, primitive directed cycles in the pointed-edge transition
graph are precisely primitive closed edge geodesics, up to cyclic shift, and
the product of their transition weights is $u^{\ell_c(P)}$.
The same statement holds for pointed chambers and closed chamber geodesics. Thus the standard determinant formula for finite weighted directed graphs (see \cite{Hashimoto1989,Bass1992}) applies to give the following result.

\begin{theorem}
\label{thm:rationality_zeta}
Let $X$ be a finite abstract affine complex of type $\widetilde A_2$. Then
\[
Z_1(X,u)=\det(\Id-A_E(u))^{-1},
\qquad
Z_2(X,u)=\det(\Id-A_C(u))^{-1}.
\]
In particular, both zeta functions are reciprocals of polynomials in
$\mathbb Z[u]$ with constant term $1$.
\end{theorem}

\section{M\"obius Vertex Polynomial}\label{s:vertex_mobius_polynomial}
In this section we construct the M\"obius vertex polynomial of $X$. We first extend the vertex set of each link $\mathrm{Lk}(x)$ to an augmented set and equip it with a partial order. The M\"obius function on this poset is then used to define a family of vertex adjacency operators, from which the M\"obius vertex polynomial is obtained by taking a determinant.

\subsection{Augmented local poset and M\"obius function}\label{s:augmented_local_poset}

Let $x\in X^{(0)}$ be a vertex of an abstract affine complex of type
$\widetilde A_2$.  
Recall from Definition~\ref{def:abstract_affine_A2} that \(\mathrm{Lk}(x)=(P_x\sqcup L_x,E_x)\) is a bipartite graph. The set $P_x$ consists of all vertices $y$ in $\mathrm{Lk}(x)$ with $\tau_d(x, y)=1$; while the set $L_x$ consists of all vertices $y$ in $\mathrm{Lk}(x)$ with $\tau_d(x, y)=2$.

\begin{definition}[Augmented local poset]
The \emph{augmented local poset} at $x$ is the set
\[
\widehat{\mathcal P}(x)
=
\{x\}
\;\cup\;
P_x
\;\cup\;
L_x
\;\cup\;
\{\pi x\}
\]
equipped with the partial order described below. Here $\pi x$ is a formal symbol representing a cyclic return to the vertex $x$;  it does not denote an actual vertex. We identify the vertex $x$ with the minimal element of $\widehat{\mathcal P}(x)$.  The partial order on
$\widehat{\mathcal P}(x)$ is the transitive closure of the relations
\(x<p<\ell<\pi x\) whenever \(p\in P_x\), \(\ell\in L_x\), and $p$
and $\ell$ are adjacent in $\mathrm{Lk}(x)$.
\end{definition}

The color difference function $\tau_d$ is extended to $\widehat{\mathcal P}(x)$ following the conventions
\(\tau_d(x,x)=0\) and \(\tau_d(x,\pi x)=3\). For clarity, whenever \(a,b\in\widehat{\mathcal P}(x)\), we write
\(a\le_x b\) to indicate that the order relation is taken in the local poset
$\widehat{\mathcal P}(x)$.
We also write $a<_x b$ for the corresponding strict order relation.









The augmented poset $\widehat{\mathcal P}(x)$ admits a natural
M\"obius function $$\mu:\widehat{\mathcal P}(x) \times \widehat{\mathcal P}(x)\to \mathbb{Z}.$$ Its restriction to the subset $\{x\}\times \widehat{\mathcal P}(x)$, viewed as a function on $\widehat{\mathcal P}(x)$, is denoted by $\mu_x$.  The function $\mu_x$ is 
characterized recursively by the relation
\[
\sum_{x\le_x y\le_x z}\mu_x(y)
=
\delta_{x,z},
\]
where $\delta$ denotes the Kronecker delta.


\subsection{The M\"obius vertex polynomial \texorpdfstring{$M(X,u)$}{M(X,u)}}\label{s:combinatorial_mobius_polynomial}

 Let \(C_0(X)\) be the free $\mathbb Z[[u]]$-module of $\mathbb Z[[u]]$-valued functions on the vertex set of $X$. Extend $f\in C_0(X)$ from $X^{(0)}$ to
\(\bigcup_{x\in X^{(0)}}\widehat{\mathcal P}(x)\) by setting
\(f(\pi x)=f(x)\).


For $i=0,1,2,3$, define the \emph{$i$--th M\"obius vertex adjacency
operator} \(\widehat A_i:C_0(X)\longrightarrow C_0(X)\) by
\[
(\widehat A_i f)(x)
=
\sum_{\substack{y\in\widehat{\mathcal P}(x)\\ \tau_d(x,y)=i}}
\mu_x(y)\,f(y).
\]

Thus $\widehat A_0=\Id$, and $\widehat A_3$ is the diagonal operator
$(\widehat A_3f)(x)=\mu_x(\pi x)f(x)$.  The operators
$\widehat A_1$ and $\widehat A_2$ locally record contributions from vertices in $P_x$ and $L_x$, respectively, in the link $\operatorname{Lk}(x)$.


We now assume that $X^{(0)}$ is finite.  Then the operators \(\widehat A_i\) for \(i=0,1,2,3\) act on the free $\mathbb Z[[u]]$-module \(C_0(X)\) of finite rank, with integer matrix entries. Let  $\widehat A(u)$ be the generating polynomial
\[
\widehat A(u)
=
\sum_{i=0}^3 \widehat A_i\,u^i
\;\in\;
\mathrm{End}_{\mathbb Z[[u]]}(C_0(X)),
\]
whose matrix has entries in $\mathbb Z[u]$.

\begin{definition}[M\"obius vertex polynomial]
We define the \emph{M\"obius vertex polynomial} of $X$ to be the
polynomial
\[
M(X,u)
=
\det\bigl(\widehat A(u)\bigr)
=
\det\!\Bigl(\sum_{i=0}^3 \widehat A_i\,u^i\Bigr).
\]
\end{definition}

Then we are ready to state the main theorem.  The proof of the theorem will be given in Section~\ref{s:pointed_operators}.

\begin{theorem}[Main zeta identity]\label{thm:main_zeta_identity}
Let $X$ be a finite abstract affine complex of type $\widetilde A_2$.
Then, in $\mathbb Z[[u]]$,
\[
\frac{Z_1(X,u)}{Z_2(X,-u)}
=
(1-u^3)^{\chi(X)}\,M(X,u)^{-1},
\]
where $\chi(X)=|X^{(0)}|-|X^{(1)}|+|X^{(2)}|$ is the ordinary Euler
characteristic of $X$ and
\[
M(X,u)
=
\det\!\Bigl(\sum_{i=0}^3\widehat A_i u^i\Bigr)
\]
is the M\"obius vertex polynomial of $X$.  Here $M(X,u)$ and
$Z_2(X,-u)$ have constant term $1$, hence are units in $\mathbb Z[[u]]$.
\end{theorem}


We end this section with a remark comparing the current framework with the classical $\mathrm{PGL}_3$ case.

\begin{remark}[The $\mathrm{PGL}_3$ quotient case]\label{rem:pgl3_mobius_values}
When $X$ arises from a finite quotient of the Bruhat--Tits building of  $\mathrm{PGL}_3$ over a local field with residue field of order $q$, each vertex link is the incidence graph of the projective plane $\mathbb{P}^2(\mathbb{F}_q)$. The augmented local poset has M\"obius values $\mu_x(x)=1$, $\mu_x(p)=-1$ for $p \in P_x$, $\mu_x(\ell)=q$ for $\ell \in L_x$, and $\mu_x(\pi x)=-q^3$.

One can easily verify these values using two standard properties of $\mathbb{P}^2(\mathbb{F}_q)$: each line is incident with $q+1$ points, and both $P_x$ and $L_x$ consist of $q^2+q+1$ elements. Specifically, for $\ell\in L_x$ the relation $\sum_{x \le_x y \le_x \ell} \mu_x(y) = 0$ yields $\mu_x(\ell) = -1 + (q+1) = q$, and the top-level sum over the entire poset yields $\mu_x(\pi x) = -q^3$.

Denote the unweighted type-one and type-two vertex adjacency operators by $A_1$ and $A_2$. Then the M\"obius vertex operator takes the explicit form
\[
\widehat A(u) = \Id - u A_1 + q u^2 A_2 - q^3 u^3 \Id.
\]
In this case, the M\"obius vertex polynomial $M(X,u) = \det\widehat A(u)$ recovers the reciprocal of the vertex $L$-function associated to $\mathrm{PGL}_3$.
\end{remark}

\section{Pointed Incidence Complex and the Proof of the Main Identity}\label{s:pointed_operators}

We now prove the main zeta identity by introducing a three-term cochain complex based on the pointed faces of $X$. Adapted from the cohomological framework of \cite{KangLiArtinPGL3}, our construction introduces the differential $d$ and auxiliary degree-lowering maps $\delta$. These maps give rise to a cochain map $\Phi$, whose degree components $\Phi_0$, $\Phi_1$, and $\Phi_2$ yield determinants that respectively encode the vertex, edge, and chamber contributions. Structurally, $\Phi$ is cochain-homotopic to multiplication by $1-u^3$, and the main zeta identity is then established via Sylvester's identity.


Throughout this section, pointed chambers $(x, y, z) $ are written in the order such that $\tau_d(x, y)=\tau_d(y, z)=\tau_d(z, x)=1$ as described in Section~\ref{s:local_geometry}.  All relations $<_{v}$ are
taken in the augmented local poset $\widehat{\mathcal P}(v)$ as described in
Section~\ref{s:augmented_local_poset}.

\subsection{The pointed cochain complex}

Let $X$ be a finite abstract affine complex of type $\widetilde A_2$. For $i=0,1,2$, let $C_i(X)$ denote the $\mathbb{Z}[[u]]$-valued  functions on vertices ($i=0$) or pointed $i$-faces ($i=1,2$). Consequently, $C_i(X)$ is a free $\mathbb{Z}[[u]]$-module of rank $(i+1)|X^{(i)}|$.

\subsubsection*{The coboundary operators $d_0$ and $d_1$}

Define the operator \(d_0:C_0(X)\longrightarrow C_1(X)\) by
\[
(d_0f)(v_0,v_1)
=
u^{\tau_d(v_0,v_1)}f(v_1)-f(v_0),
\]
and the operator \(d_1:C_1(X)\longrightarrow C_2(X)\) by
\[
(d_1f)(v_0,v_1,v_2)
=
u\,f(v_1,v_2)-f(v_0,v_2)+f(v_0,v_1).
\]

These are weighted incidence maps. The power of $u$ records the color increments of the pointed faces. A direct computation gives $d_1d_0=0$.  Indeed,
\begin{align*}
(d_1d_0f)(v_0,v_1,v_2)
&=
u(d_0f)(v_1,v_2)
-(d_0f)(v_0,v_2)
+(d_0f)(v_0,v_1)\\
&=
u\bigl(uf(v_2)-f(v_1)\bigr)
-\bigl(u^2f(v_2)-f(v_0)\bigr)\\
&\quad
+\bigl(uf(v_1)-f(v_0)\bigr)\\
&=0.
\end{align*}
Thus
\[
0\longrightarrow C_0(X)
\xrightarrow{d_0}
C_1(X)
\xrightarrow{d_1}
C_2(X)
\longrightarrow 0
\]
is a cochain complex.  

\subsubsection*{The auxiliary operators \texorpdfstring{$\delta_1$}{delta 1} and \texorpdfstring{$\delta_2$}{delta 2}}

We also need two auxiliary operators $\delta_1: C_1(X)\longrightarrow C_0(X)$ and $\delta_2:C_2(X)\longrightarrow C_1(X)$ defined as follows.

For $f\in C_1(X)$ and $v_0\in X^{(0)}$, define
\[
(\delta_1f)(v_0)
=
-\sum_{v_0<_{{v_0}}x<_{{v_0}}\pi v_0}
u^{\tau_d(v_0,x)}
\mu_{v_0}(x)\,
f(x,v_0).
\]

For $f\in C_2(X)$, if $(v_0,v_2)$ is a pointed edge of type $2$ so that
$v_2\in L_{v_0}$, define
\[
(\delta_2f)(v_0,v_2)
=
-u\sum_{v_0<_{{v_0}}x<_{{v_0}}v_2}
f(x,v_2,v_0); 
\]

while if $(v_0,v_1)$ is a pointed edge of type $1$ so that
$v_1\in P_{v_0}$, define
\[
(\delta_2f)(v_0,v_1)
=
u^2\sum_{v_1<_{{v_1}}x<_{{v_1}}v_0}
f(x,v_0,v_1).
\]

The two cases differ because of the ordering convention for pointed
chambers.  With distinguished vertex $x$ opposite the edge, the chamber
is $(x,v_0,v_1)$ for a type $1$ edge but $(x,v_2,v_0)$ for a type $2$
edge, in which the edge appears reversed.  This reversal accounts for
the minus sign and the different power of $u$.

Observe that if $(v_0,v_2)$ has type $2$, then $(v_2, v_0)$ has type $1$, and it follows from the above definition that 
\((\delta_2f)(v_2,v_0)=-u(\delta_2f)(v_0,v_2)\).
This shows the compatibility of the two edge-type formulas
under reversal.

\subsubsection*{The cochain map $\Phi$}

Define the homomorphisms $\Phi_i: C_i(X) \longrightarrow C_i(X)$ as follows:
\[
\begin{aligned}
\Phi_0&=(1-u^3)\Id_{C_0(X)}+\delta_1d_0,\\
\Phi_1&=(1-u^3)\Id_{C_1(X)}+d_0\delta_1+\delta_2d_1,\\
\Phi_2&=(1-u^3)\Id_{C_2(X)}+d_1\delta_2.
\end{aligned}
\]

On $C_0(X)\oplus C_1(X)\oplus C_2(X)$, write $d=d_0+d_1$, $\delta=\delta_1+\delta_2$, and let $\Phi$ denote the
operator with degree components $\Phi_0,\Phi_1,\Phi_2$.  Then$$
\Phi=(1-u^3)\Id+d\delta+\delta d.$$
It can be easily verified that $\Phi_1 d_0=d_0\Phi_0$ and $\Phi_2 d_1=d_1 \Phi_1$. So $\Phi$ is a cochain map that is cochain-homotopic to the scalar operator $(1-u^3)\Id$.  

We next study the three determinants $\det(\Phi_0), \det(\Phi_1)$, and $\det(\Phi_2)$, and show that they are related to vertex, edge, and chamber zeta functions, respectively.





\subsection{Computation of \texorpdfstring{$\Phi_0$}{Phi 0}}\label{s:Phi0_computation}

The degree-zero calculation is the point at which the present argument
differs most visibly from the regular $\mathrm{PGL}_3$ case: the
uniform coefficients are replaced by the M\"obius functions of the
augmented local posets.  Nevertheless, the computation is immediate
from the defining M\"obius relation.  For $f\in C_0(X)$ and
$v_0\in X^{(0)}$,
\[
\begin{aligned}
(\Phi_0f)(v_0)
&=(1-u^3)f(v_0)
-\sum_{v_0<_{v_0}x<_{v_0}\pi v_0}
 u^{\tau_d(v_0,x)}\mu_{v_0}(x)
 \bigl(u^{\tau_d(x,v_0)}f(v_0)-f(x)\bigr)\\
&=f(v_0)
+\sum_{v_0<_{v_0}x<_{v_0}\pi v_0}
 u^{\tau_d(v_0,x)}\mu_{v_0}(x)f(x)
+u^3\mu_{v_0}(\pi v_0)f(v_0).
\end{aligned}
\]
Here we used \(\tau_d(v_0,x)+\tau_d(x,v_0)=3\) and
\(\sum_{v_0\le_{v_0}x\le_{v_0}\pi v_0}\mu_{v_0}(x)=0\).
By the definition of the M\"obius vertex adjacency operators, this is
precisely
\[
\Phi_0=\sum_{i=0}^3\widehat A_i u^i=\widehat A(u).
\]
Consequently,
\[
\boxed{\det(\Phi_0)=M(X,u).}
\]

\subsection{Computation of \texorpdfstring{$\Phi_2$}{Phi 2}}\label{s:Phi2_computation}

The chamber calculation has the same formal structure as the degree-two
calculation in~\cite[Sections~5--6]{KangLiArtinPGL3}.  We record the
factorization in the present notation and indicate the local
coefficient check, since it also fixes the sign in $Z_2(X,-u)$.

Define the chamber rotation operator by
\((J_Cf)(v_0,v_1,v_2)=f(v_1,v_2,v_0)\).
Then $J_C^3=\Id$.  If $c=(v_0,v_1,v_2)$, the panel-successor convention
gives
\[
(J_CA_C(u)J_C^{-1}f)(c)
=
u\sum_{\substack{v_0<_{v_0}x<_{v_0}v_2\\x\neq v_1}}
f(x,v_2,v_0).
\]

\begin{proposition}\label{prop:Phi2_factorization}
On $C_2(X)$,
\[
\Phi_2
=
(\Id+uJ_C+u^2J_C^2)
(\Id+J_CA_C(u)J_C^{-1}).
\]
\end{proposition}

\begin{proof}
For a pointed chamber $c=(v_0,v_1,v_2)$, direct expansion of
$d_1\delta_2f$ gives the three full edge-residue sums
\[
\begin{aligned}
(d_1\delta_2f)(c)
&=
u\sum_{v_0<_{v_0}x<_{v_0}v_2}f(x,v_2,v_0)\\
&\quad+
u^2\sum_{v_1<_{v_1}x<_{v_1}v_0}f(x,v_0,v_1)\\
&\quad+
u^3\sum_{v_2<_{v_2}x<_{v_2}v_1}f(x,v_1,v_2).
\end{aligned}
\]
Expanding the proposed product gives these same sums with, respectively,
the terms $v_1,v_2,v_0$ omitted, together with
\[
f(c)+u f(J_Cc)+u^2f(J_C^2c).
\]
After $(1-u^3)f(c)$ is incorporated in
$\Phi_2=d_1\delta_2+(1-u^3)\Id$, these three displayed terms supply
exactly the omitted summands.  This proves the factorization.
\end{proof}

Let $N_2=|X^{(2)}|$.  The pointed chambers decompose into $N_2$
rotation orbits of the form \(\{c,J_Cc,J_C^2c\}\).
On each orbit, $\Id+uJ_C+u^2J_C^2$ has matrix
\[
\begin{pmatrix}
1&u&u^2\\
u^2&1&u\\
u&u^2&1
\end{pmatrix},
\]
whose determinant is $(1-u^3)^2$.  Moreover,
\(\det(\Id+J_CA_C(u)J_C^{-1})=\det(\Id+A_C(u))\), and
$A_C(-u)=-A_C(u)$ because every chamber successor step has weight
$u$.  Proposition~\ref{prop:Phi2_factorization} and
Theorem~\ref{thm:rationality_zeta} therefore give
\[
\boxed{
\det(\Phi_2)
=
(1-u^3)^{2N_2}\det(\Id+A_C(u))
=
(1-u^3)^{2N_2}Z_2(X,-u)^{-1}.
}
\]

\subsection{Computation of \texorpdfstring{$\Phi_1$}{Phi 1}}\label{s:Phi1_computation}

The degree-one calculation is the longest local coefficient check in the
determinant argument.  We follow the general approach of
\cite[Sections~5--6]{KangLiArtinPGL3}, adapting the local coefficient computation to the present setting through the M\"obius relation and the weak $A_2$ link axioms. 

Our strategy is to reorganize the calculation according to the two
types of pointed edges. We first introduce auxiliary operators that
separate the contributions of the two edge types and account for
the reversal of pointed edges. The resulting operator identities
will then reduce the computation of $\Phi_1$ to a determinant
associated with the edge adjacency operator.

Write \(C_1(X)=C_1^{(1)}(X)\oplus C_1^{(2)}(X)\) according to
pointed-edge type.  Define the unipotent operator $Q$ by
\((Qf)(v_0,v_1)=f(v_0,v_1)\) on type-$1$ pointed edges and
\[
(Qf)(v_0,v_2)
=
f(v_0,v_2)
+
u\sum_{v_0<_{v_0}x<_{v_0}v_2}f(x,v_2)
\]
on type-$2$ pointed edges.  Relative to the type decomposition,
\[
Q=
\begin{pmatrix}
\Id&0\\
Q_{21}&\Id
\end{pmatrix},
\qquad
\det Q=1.
\]

Define weighted edge reversal by
\[
(J_Ef)(v_0,v_1)
=
u^{\tau_d(v_0,v_1)}f(v_1,v_0),
\qquad
W=\Id+J_E.
\]
Then $J_E^2=u^3\Id$, so $J_E$ is not invertible over $\mathbb Z[[u]]$.
Instead, let $\widetilde A_E(u)$ be the operator on $C_1(X)$ given by
\[
(\widetilde A_E(u)f)(v_0,v_1)
=
u^2\sum_{\substack{x\in L_{v_0}\\x\text{ not incident with }v_1}}
f(x,v_0)
\]
for a type-$1$ pointed edge, and
\[
(\widetilde A_E(u)f)(v_0,v_2)
=
u\sum_{\substack{x\in P_{v_0}\\x\text{ not incident with }v_2}}
f(x,v_0)
\]
for a type-$2$ pointed edge.  Lemma~\ref{lem:spherical_A2_basic}(2) and edge reversal give
\[
\widetilde A_E(u)J_E=J_EA_E(u).
\]

Finally, define the residual operator
\[
(\mathcal Nf)(v_0,v_1)
=
\begin{cases}
\displaystyle
u\sum_{\substack{x\in P_{v_0}\\x\neq v_1}}f(x,v_0),
&\tau_d(v_0,v_1)=1,\\[2ex]
0,
&\tau_d(v_0,v_1)=2.
\end{cases}
\]

\begin{proposition}[Degree-one local operator identities]
\label{prop:degree_one_local_identities}
For $f\in C_1(X)$, $g\in C_0(X)$, and $h\in C_2(X)$,
\[
(\delta_1Qf)(v_0)
=
\sum_{v_0<_{v_0}x<_{v_0}\pi v_0}
u^{\tau_d(v_0,x)}f(x,v_0),
\]
\[
(Wd_0g)(v_0,v_1)=-(1-u^3)g(v_0),
\]
and
\[
(W\delta_2h)(v_0,v_i)
=
\begin{cases}
0,&i=1,\\
(1-u^3)(\delta_2h)(v_0,v_2),&i=2.
\end{cases}
\]
These identities imply
\begin{equation}\label{eq:Phi1_conjugation_identity}
\frac{1}{1-u^3}W\Phi_1Q
=
\Id-\widetilde A_E(u)-\mathcal N.
\end{equation}
\end{proposition}

The proof of this proposition relies on two local properties of the
links. The first identifies the relevant vertices in the links at
the two endpoints of a pointed edge, allowing the local sums arising
in the computation to be expressed in a common link. The second
is the unique common-neighbor property in a weak $A_2$ link, which
will be used to reorganize the double sum arising from pointed edges
of type $2$. We record these two facts before proving the proposition.


\begin{lemma}[Edge-residue identification]
\label{lem:edge_residue_identification}
Let $(v_0,v_1)$ be a pointed edge of type $1$.
Then the following sets are naturally in bijection:
\begin{enumerate}
    \item the vertices $x\in L_{v_0}$ adjacent to $v_1\in P_{v_0}$;
    \item the vertices $x\in P_{v_1}$ adjacent to $v_0\in L_{v_1}$.
\end{enumerate}
Equivalently, the chambers containing the edge $\{v_0,v_1\}$ may be
parametrized either from the link at $v_0$ or from the link at $v_1$.

More precisely, if
\[
    v_0<_{v_0}v_1<_{v_0}x<_{v_0}\pi v_0,
\]
then the same chamber may be written from the link at $v_1$ as
\[
    v_1<_{v_1}x<_{v_1}v_0<_{v_1}\pi v_1.
\]
The corresponding statement for pointed edges of type $2$ is obtained
by reversing the pointed edge.
\end{lemma}

\begin{proof}
Both sets in the statement are canonically identified with the set of
chambers containing the edge $\{v_0,v_1\}$. Indeed, a vertex
$x\in L_{v_0}$ adjacent to $v_1\in P_{v_0}$ means precisely that
$\{v_0,v_1,x\}$ is a chamber. Viewed from $v_1$, the same vertex $x$
belongs to $P_{v_1}$, while $v_0$ belongs to $L_{v_1}$; moreover,
$x$ is adjacent to $v_0$ in $\mathrm{Lk}(v_1)$. This gives the stated
identification.
\end{proof}

\begin{lemma}[Local $A_2$ adjacency]\label{lem:local_projective_plane_incidence}
Let $x\in X^{(0)}$. In the bipartite graph
\[
    \mathrm{Lk}(x)=(P_x\sqcup L_x,E_x),
\]
any two distinct vertices in $P_x$ have a unique common neighbor in
$L_x$, and any two distinct vertices in $L_x$ have a unique common
neighbor in $P_x$.

Consequently, for any $\ell\in L_x$, the map
\[
\begin{aligned}
\big\{(p,m)\in P_x\times L_x:\;&
p\sim\ell,\ p\sim m,\ m\ne\ell\big\}
&\longrightarrow L_x\setminus\{\ell\},\\
(p,m)&\longmapsto m
\end{aligned}
\]
is a bijection. The dual statement, with $P_x$ and $L_x$
interchanged, also holds.
\end{lemma}

\begin{proof}
The first assertion follows from Lemma~\ref{lem:spherical_A2_basic} applied to the weak
$A_2$ link $\mathrm{Lk}(x)$. For the bijection, fix $m\in L_x\setminus\{\ell\}$. By the uniqueness assertion, $\ell$ and $m$ have a unique common neighbor
$p\in P_x$. Thus $p\sim\ell$ and $p\sim m$, so $(p,m)$ is the unique
preimage of $m$. This proves the bijection. The dual statement follows
in the same way.
\end{proof}

With these two local facts in hand, we return to the proof of Proposition~\ref{prop:degree_one_local_identities}.

\begin{proof}
For the first identity, expand $\delta_1Q$ and collect the coefficient
of each local point or line.  A local point $p$ has
$\mu_{v_0}(p)=-1$, so its coefficient is $u$.  For a local line
$\ell$, the coefficient is
\[
-u^2\sum_{v_0<_{v_0}x\le_{v_0}\ell}\mu_{v_0}(x)=u^2
\]
by the defining M\"obius relation.  The second identity follows by
expanding $Wd_0$ and using
\[
\tau_d(v_0,v_1)+\tau_d(v_1,v_0)=3.
\]
For a type-$1$ pointed edge, the two formulas defining $\delta_2$ give
\[
(\delta_2h)(v_0,v_1)+u(\delta_2h)(v_1,v_0)=0.
\]
For a type-$2$ pointed edge, they give
\[
(\delta_2h)(v_0,v_2)+u^2(\delta_2h)(v_2,v_0)
=(1-u^3)(\delta_2h)(v_0,v_2).
\]
These are precisely the two cases of the third identity after weighted
edge reversal.

It remains to identify the surviving terms in $W\Phi_1Q$.  On a
type-$1$ pointed edge, the preceding identities give
\[
\begin{aligned}
\frac{1}{1-u^3}(W\Phi_1Qf)(v_0,v_1)
&=
f(v_0,v_1)
-u^2\sum_{\substack{x\in L_{v_0}\\x\text{ not incident with }v_1}}
 f(x,v_0)\\
&\quad
-u\sum_{\substack{x\in P_{v_0}\\x\neq v_1}}f(x,v_0).
\end{aligned}
\]
Here the terms indexed by local lines in $L_{v_0}$ incident with
$v_1$ cancel against the corresponding terms indexed by local points
in $P_{v_1}$ incident with $v_0$ in $WQf$ and $\delta_1Qf$.
Lemma~\ref{lem:edge_residue_identification} identifies these two
index sets: both parametrize the chambers containing the edge
$\{v_0,v_1\}$. Thus the cancellation is a direct application of the
edge-residue identification.

On a type-$2$ pointed edge $(v_0,v_2)$, apply
Lemma~\ref{lem:edge_residue_identification} to the reversed
type-$1$ pointed edge $(v_2,v_0)$. It identifies the local lines
$x\in L_{v_2}$ incident with $v_0$ with the local points
$x\in P_{v_0}$ incident with $v_2$, so the corresponding chamber
contributions can be reindexed in the link at $v_0$. Now
Lemma~\ref{lem:local_projective_plane_incidence} gives a bijection
between the local lines $m\in L_{v_0}\setminus\{v_2\}$ and the pairs
$(p,m)$ with $p\in P_{v_0}$ incident with both $v_2$ and $m$.
Consequently, each local line contribution is paired with exactly
one local point contribution, and these terms cancel. This leaves
\[
\frac{1}{1-u^3}(W\Phi_1Qf)(v_0,v_2)
=
f(v_0,v_2)
-u\sum_{\substack{x\in P_{v_0}\\x\text{ not incident with }v_2}}
 f(x,v_0).
\]
The two displayed formulas are exactly
\eqref{eq:Phi1_conjugation_identity} in the two edge types.
\end{proof}

We now take determinants.  Multiplying
\eqref{eq:Phi1_conjugation_identity} on the left by $\Id-J_E$ and using
\[
(\Id-J_E)(\Id+J_E)=(1-u^3)\Id
\]
gives
\[
\Phi_1
=
(\Id-J_E)
(\Id-\widetilde A_E(u)-\mathcal N)Q^{-1}.
\]
The operator $\widetilde A_E(u)$ preserves edge type.  Denote its
restrictions to $C_1^{(1)}(X)$ and $C_1^{(2)}(X)$ by $u^2B_1$ and
$uB_2$, respectively.  The operator $\mathcal N$ maps the type-$2$
summand into the type-$1$ summand; write $\mathcal N_{12}$ for its only
possibly nonzero block.  Consequently,
\[
\widetilde A_E(u)+\mathcal N
=
\begin{pmatrix}
B_1u^2&\mathcal N_{12}\\
0&B_2u
\end{pmatrix},
\]
and hence
\[
\det(\Id-\widetilde A_E(u)-\mathcal N)
=
\det(\Id-A_E(u)).
\]
Indeed, the off-diagonal block $\mathcal N_{12}$ does not affect the
determinant, and taking determinants in
$(\Id-\widetilde A_E(u))J_E=J_E(\Id-A_E(u))$ gives
$\det(\Id-\widetilde A_E(u))\det J_E=\det J_E\det(\Id-A_E(u))$, where
$\det J_E=(-u^3)^{|X^{(1)}|}$ is a nonzero element of the integral
domain $\mathbb Z[[u]]$ and can be cancelled.
Each unpointed edge contributes to $\Id-J_E$ the block
\[
\begin{pmatrix}
1&-u^{\tau_d(v_1,v_0)}\\
-u^{\tau_d(v_0,v_1)}&1
\end{pmatrix},
\]
whose determinant is $1-u^3$.  Therefore, with
$N_1=|X^{(1)}|$,
\[
\boxed{
\det(\Phi_1)
=
(1-u^3)^{N_1}\det(\Id-A_E(u))
=
(1-u^3)^{N_1}Z_1(X,u)^{-1}.
}
\]
\subsection{Proof of the main zeta identity}

Write $N_i=|X^{(i)}|$ and put $s=1-u^3$.  Since each unpointed
edge has two pointings and each chamber has three, the pointed modules have
ranks
\[
\operatorname{rank} C_0=N_0,
\qquad
\operatorname{rank} C_1=2N_1,
\qquad
\operatorname{rank} C_2=3N_2.
\]

It remains to combine the three determinant computations.  On $C_1(X)$
set
\[
A=d_0\delta_1,
\qquad
B=\delta_2d_1.
\]
Since $d_1d_0=0$,
\[
BA=\delta_2d_1d_0\delta_1=0.
\]
Hence
\[
s^{-1}\Phi_1
=
(\Id+s^{-1}B)(\Id+s^{-1}A).
\]
Taking determinants and applying Sylvester's identity to the two adjacent
incidence maps gives
\[
\begin{aligned}
s^{-2N_1}\det(\Phi_1)
&=
\det(\Id+s^{-1}\delta_2d_1)
\det(\Id+s^{-1}d_0\delta_1)\\
&=
\det(\Id+s^{-1}d_1\delta_2)
\det(\Id+s^{-1}\delta_1d_0)\\
&=
s^{-3N_2-N_0}\det(\Phi_2)\det(\Phi_0).
\end{aligned}
\]
Thus
\[
\det(\Phi_0)\det(\Phi_2)
=
s^{N_0-2N_1+3N_2}\det(\Phi_1).
\]
This is the determinant manifestation, in our three-term complex, of the
cochain homotopy
$\Phi\simeq s\Id$.

The computations in Sections~\ref{s:Phi0_computation},
\ref{s:Phi2_computation}, and~\ref{s:Phi1_computation} give
\[
\det(\Phi_0)=M(X,u),
\]
\[
\det(\Phi_1)=s^{N_1}Z_1(X,u)^{-1},
\qquad
\det(\Phi_2)=s^{2N_2}Z_2(X,-u)^{-1}.
\]
Substitution yields
\[
M(X,u)s^{2N_2}Z_2(X,-u)^{-1}
=
s^{N_0-N_1+3N_2}Z_1(X,u)^{-1}.
\]
Hence
\[
\frac{Z_1(X,u)}{Z_2(X,-u)}
=
s^{N_0-N_1+N_2}M(X,u)^{-1}.
\]
Since $N_0-N_1+N_2=\chi(X)$ and $s=1-u^3$, we obtain
\[
\boxed{
\frac{Z_1(X,u)}{Z_2(X,-u)}
=
(1-u^3)^{\chi(X)}M(X,u)^{-1}.
}
\]
This proves Theorem~\ref{thm:main_zeta_identity}.

\section{Zeta identity verification for the examples}\label{s:example_verifications}

We close with determinant checks for the examples.  These computations
are meant to confirm the normalizations in concrete finite cases, not to
give separate proofs of the main theorem.  The displayed factorizations
are obtained directly from the finite vertex, edge successor, and
chamber successor matrices attached to the corresponding incidence data.

For each fixed complex $X$, write
\[
\begin{aligned}
M(X,u)&=\det\widehat A(u),\\
D_E(u)&=\det(\Id-A_E(u)),\\
D_C^+(u)&=\det(\Id+A_C(u))=\det(\Id-A_C(-u)).
\end{aligned}
\]
All four examples below have Euler characteristic zero.  Thus the zeta
identity in each case is equivalent to the finite determinant relation
\(D_E(u)=M(X,u)D_C^+(u)\).

\begin{example}[Determinant check for the minimal apartment torus]
Let $X_{\mathrm{tor}}$ be the minimal apartment torus from
Example~\ref{ex:minimal-apartment-torus} and put $t=u^3$.  Then
\(M(X_{\mathrm{tor}},u)=\det \widehat A(u)=(1-t)^9\).
\[
D_E(u)
=(t-1)^{18}(t+1)^9
=(1-t)^{18}(1+t)^9.
\]
The chamber successor permutation consists of nine cycles of length $6$,
so \(D_C^+(u)=(1-t^2)^9\).
Consequently,
\[
D_E(u)
=
(1-t)^9(1-t^2)^9
=
M(X_{\mathrm{tor}},u)D_C^+(u).
\]
\end{example}

\begin{example}[The sign in the Klein bottle quotient]
Let $X=X_{\mathrm{Kl}}$ be the Klein bottle quotient from
Example~\ref{ex:klein-bottle-quotient}.  Its chamber successor operator
has the form $A_C(u)=uB_C$, where $B_C$
is a permutation matrix.  The successor permutation on pointed chambers
has cycle lengths \(9,9,18,72\).
Therefore
\[
\det(\Id-A_C(u))
=
(1-u^9)^2(1-u^{18})(1-u^{72}),
\]
whereas
\[
D_C^+(u)=\det(\Id+A_C(u))=\det(\Id-A_C(-u))
=
(1+u^9)^2(1-u^{18})(1-u^{72}).
\]
Thus $Z_2(X,u)$ is not invariant under $u\mapsto -u$; indeed
\[
\frac{Z_2(X,-u)}{Z_2(X,u)}
=
\left(\frac{1-u^9}{1+u^9}\right)^2.
\]
This explains why the main identity is naturally normalized with
$Z_2(X,-u)$.

For this quotient, the vertex determinant is
\(M(X,u)=(1-u^9)^4(1+u^{18})\).  The directed edge successor cycles
have total weights \(u^9,u^9,u^{18},u^{18},u^{36},u^{72}\).
Hence
\[
D_E(u)
=
(1-u^9)^2(1-u^{18})^2(1-u^{36})(1-u^{72}).
\]
A direct comparison gives \(D_E(u)=M(X,u)D_C^+(u)\).
\end{example}

\begin{example}[Determinant check for the glued two-torus example]
Let $X$ be the two-torus gluing example from
Example~\ref{ex:two-torus-gluing}, and put $t=u^3$.  A direct calculation
gives
\[
M(X,u)
=
(t-1)^5(t^2+1)(8t^2-5t+1)(8t^3-11t^2+6t-1).
\]
Similarly,
\[
\begin{aligned}
D_E(u)
={}&
(t-1)^{12}(t+1)^7(t^2+1)^2(t^4+1) (8t^2-5t+1)(8t^4-5t^2+1) \\
&(8t^3-11t^2+6t-1) (8t^6-11t^4+6t^2-1),
\end{aligned}
\]
\[
\begin{aligned}
D_C^+(u)
={}&
(t-1)^7(t+1)^7(t^2+1)(t^4+1)  (8t^4-5t^2+1)(8t^6-11t^4+6t^2-1).
\end{aligned}
\]
These factors satisfy \(D_E(u)=M(X,u)D_C^+(u)\).
\end{example}

\begin{example}[Determinant check for the torus--Klein gluing example]
Let $X$ be the torus--Klein gluing example from
Example~\ref{ex:torus-klein-gluing}.  Its two link types give the local
M\"obius values $\mu_x(\pi x)=-1$ at the hexagonal-link vertices and
$\mu_x(\pi x)=-2$ at the vertices on the common strip.
The following finite computation checks the determinant normalization for
this example.  The matrices $\widehat A(u)$, $A_E(u)$, and $A_C(u)$ are
formed, respectively, from the local M\"obius data of
Section~\ref{s:vertex_mobius_polynomial} and the directed edge and chamber
successor relations.
Their sizes are \(27\), \(2|X^{(1)}|=198\), and
\(3|X^{(2)}|=216\), respectively.  These matrices are determined
entirely by the incidence
data of Example~\ref{ex:torus-klein-gluing}; in particular, they use the
inherited chamber set and no flag-completion chambers are added.  A
direct calculation gives \(D_E(u)=M(X,u)D_C^+(u)\).
The explicit determinant factors are omitted because of their size.
Thus
\[
\frac{Z_1(X,u)}{Z_2(X,-u)}
=
\frac{D_C^+(u)}{D_E(u)}
=
M(X,u)^{-1}
=
(1-u^3)^{\chi(X)}M(X,u)^{-1}.
\]
\end{example}

\section*{Declaration on the Use of AI}
The authors used ChatGPT and Claude to assist with reviewing the manuscript and refining the text. The authors critically reviewed and verified all AI-assisted suggestions and take full responsibility for the content of the published article.

\begingroup
\sloppy
\hbadness=10000
\bibliographystyle{plain}
\bibliography{references}

@article{Bass1992,
  author = {Bass, Hyman},
  title = {The {Ihara}-{Selberg} zeta function of a tree lattice},
  journal = {Internat. J. Math.},
  volume = {3},
  number = {6},
  year = {1992},
  pages = {717--797}
}

@article{DeitmarKang2013,
  author = {Deitmar, Anton and Kang, Ming-Hsuan},
  title = {Geometric zeta functions for higher rank {$p$}-adic groups},
  journal = {Illinois J. Math.},
  volume = {58},
  number = {3},
  year = {2014},
  pages = {719--738}
}

@article{FangLiWang2013,
  author = {Fang, Yang and Li, Wen-Ching Winnie and Wang, Chian-Jen},
  title = {The zeta functions of complexes from {$Sp(4)$}},
  journal = {Int. Math. Res. Not. IMRN},
  volume = {2013},
  number = {4},
  year = {2013},
  pages = {886--923}
}

@incollection{Hashimoto1989,
  author = {Hashimoto, Ki-ichiro},
  title = {Zeta functions of finite graphs and representations of {$p$}-adic groups},
  booktitle = {Automorphic Forms and Geometry of Arithmetic Varieties},
  series = {Adv. Stud. Pure Math.},
  volume = {15},
  publisher = {Academic Press},
  address = {Boston, MA},
  year = {1989},
  pages = {211--280}
}

@misc{HongKwon2025,
  author = {Hong, Soonki and Kwon, Sanghoon},
  title = {Chamber zeta function and closed galleries in the standard non-uniform complex from {$\mathrm{PGL}_3$}},
  year = {2025},
  note = {arXiv:2512.23276}
}

@article{Ihara1966,
  author = {Ihara, Yasutaka},
  title = {On discrete subgroups of the two by two projective linear group over {$p$}-adic fields},
  journal = {J. Math. Soc. Japan},
  volume = {18},
  year = {1966},
  pages = {219--235}
}

@article{KangLi2014,
  author = {Kang, Ming-Hsuan and Li, Wen-Ching Winnie},
  title = {Zeta functions of complexes arising from {$\mathrm{PGL}_3$}},
  journal = {Adv. Math.},
  volume = {256},
  year = {2014},
  pages = {46--103}
}

@article{KangLiArtinPGL3,
  author = {Kang, Ming-Hsuan and Li, Wen-Ching Winnie},
  title = {Artin {$L$}-Functions on Finite Quotients of {$\mathrm{PGL}_3$}},
  journal = {International Mathematics Research Notices},
  volume = {2015},
  number = {19},
  year = {2015},
  pages = {9251--9276},
  doi = {10.1093/imrn/rnu215}
}

@article{KangLiWang2010,
  author = {Kang, Ming-Hsuan and Li, Wen-Ching Winnie and Wang, Chian-Jen},
  title = {The zeta functions of complexes from {$\mathrm{PGL}_3$}: a representation-theoretic approach},
  journal = {Israel J. Math.},
  volume = {177},
  year = {2010},
  pages = {335--348}
}

@article{KangLiWang2018,
  author = {Kang, Ming-Hsuan and Li, Wen-Ching Winnie and Wang, Chian-Jen},
  title = {Zeta and {$L$}-functions of finite quotients of apartments and buildings},
  journal = {Israel J. Math.},
  volume = {228},
  number = {1},
  year = {2018},
  pages = {79--117},
  doi = {10.1007/s11856-018-1756-3}
}

@misc{KangYu2026,
  author = {Kang, Ming-Hsuan and Yu, Jiu-Kang},
  title = {Zeta functions of {$\mathrm{PGL}_n$} over non-{Archimedean} local fields},
  year = {2026},
  note = {\href{https://arxiv.org/abs/2607.21262}{arXiv:2607.21262 [math.NT]}}
}

@article{LLP2020,
  author = {Lubetzky, Eyal and Lubotzky, Alexander and Parzanchevski, Ori},
  title = {Random walks on {Ramanujan} complexes and digraphs},
  journal = {J. Eur. Math. Soc.},
  volume = {22},
  number = {11},
  year = {2020},
  pages = {3441--3466},
  doi = {10.4171/JEMS/990}
}

@article{LSV2005Explicit,
  author = {Lubotzky, Alexander and Samuels, Beth and Vishne, Uzi},
  title = {Explicit constructions of {Ramanujan} complexes of type {$\widetilde A_d$}},
  journal = {European J. Combin.},
  volume = {26},
  number = {6},
  year = {2005},
  pages = {965--993}
}

@article{LSV2005Ramanujan,
  author = {Lubotzky, Alexander and Samuels, Beth and Vishne, Uzi},
  title = {Ramanujan complexes of type {$\widetilde A_d$}},
  journal = {Israel J. Math.},
  volume = {149},
  year = {2005},
  pages = {267--299}
}
\endgroup

\end{document}